\documentclass{birkau}

\usepackage[utf8]{inputenc} 
\usepackage{amsmath,amssymb,url}
\usepackage{hyperref}
\usepackage{amsthm}
\usepackage{graphicx}
\usepackage{tikz}
\usepackage{tikz-cd}
\usepackage{stmaryrd}
\usepackage{mathrsfs,upgreek}
\usepackage{eucal}
\usepackage{enumerate}
\usepackage{changes}
\usepackage{verbatim}
\usepackage[capitalise]{cleveref}

\numberwithin{equation}{section}

\newtheorem{Thm}{Theorem}[section]

\newtheorem{Lem}[Thm]{Lemma}
\newtheorem{Cor}[Thm]{Corollary}
\theoremstyle{definition}

\newtheorem{Ex}[Thm]{Example}
\theoremstyle{definition}
\newtheorem{Def}[Thm]{Definition}
\numberwithin{equation}{section}

\newcommand{\op}[1]{\textrm{\upshape #1}}

\newcommand{\join}{\vee}

\newcommand{\meet}{\wedge}

\newcommand{\alg}[1]{{\textbf{\upshape #1}}}  %
\newcommand{\vv}[1]{\mathcal {#1}}

\renewcommand{\a}{\alpha}
\renewcommand{\b}{\beta}
\renewcommand{\d}{\delta}

\newcommand{\g}{\gamma}

\renewcommand{\th}{\theta}

\newcommand{\NN}{\mathbb{N}}

\newcommand{\ib}{\item[$\bullet$]}

\newcommand{\Con}{\operatorname{\mathcal{CON}}}
\newcommand{\con}{\operatorname{Con}}
\newcommand{\Cg}{\operatorname{Cg}}
\newcommand{\Eq}{\operatorname{Eq}}
\newcommand{\adm}{\operatorname{\mathcal{A}
}}

\newcommand{\vuc}[2]{#1_1,\dots,#1_{#2}}

\begin{document}
	\markboth{Stefano Fioravanti}
	{On some admissible lattices}
	
	%
	%
	
	\title{On some admissible lattices}

    \author[S.~Fioravanti]{Stefano Fioravanti}
\address{Department of Algebra\\ Charles University\\Praha\\ Czechia
}
\email{\tt stefano.fioravanti66@gmail.com}
	
\thanks{Supported by Charles University Prague under grant number PRIMUS/24/SCI/008}

	\begin{abstract} This paper explores applications of the so-called \emph{Freese's technique}, a classical approach to study the congruence variety of a given algebra. We leverage this tool to investigate lattices that are admissible in a given congruence variety. In particular, we present characterizations of congruence modular varieties, Taylor varieties, and varieties satisfying a non-trivial congruence identity by means of lattice omission provable using Freese's technique.
	\end{abstract}

    	\keywords{Congruence lattices, sublattices, congruence varieties, Freese's technique}
	
	\subjclass{06B10,  06B99}

	\maketitle

	\section{Introduction}
	
	The motivation for this work arises from our interest in a result by R. Freese and P. Lipparini. In the proof of this result, stated as Theorem \ref{freese}, the authors consider the congruence lattice of a subpower of a given algebra whose underlying set is defined by tuples of elements in a fixed congruence of the original algebra. This approach, which was first denoted as \emph{Freese's technique} in \cite{ABF}, serves as a powerful tool to identify specific sublattices within the congruence lattices of algebras in a variety. This type of problem can be considered classical in general algebra and originates from the 
    characterizations of modular lattices \cite{Dedekind1900} and distributive lattices \cite{Birkhoff1935OnTS}. Since then, researchers have studied the congruence lattice of algebras in terms of lattice omissions; see \cite{Jnsson1968} or \cite{Freese1995}, where the authors proved that congruence lattices and lattices are substantially different objects in the sense that there are varieties of lattices that are not congruence varieties.

    Recently, the omission of lattices was studied in \cite{Freese2024}, where the authors proved that if a variety satisfies a non-trivial congruence identity, then its congruence variety is finitely based only if it is distributive. Furthermore, recent attempts have also used computer-based approaches, mainly using neural networks, with the aim of learning to relate equational properties of the lattices and the omission of some lattices \cite{Ai4ua2, ai4ua1}. 
	
	Inspired by the methodology used in \cite{Freese1995} and \cite{Freese2024}, we continue the investigation initiated in \cite{ABF}, studying lattice omission in the context of congruence modular varieties, Taylor varieties, and varieties satisfying a non-trivial congruence identity.

    The main goal of this paper is to further develop the theory of applications of Freese's technique by applying it to several of the most significant lattices in the literature used to characterize Mal'cev conditions. The first of our main results regards the presence of particular sublattices in the congruence variety generated by a non-modular variety. We find an analogue of  \cite{Freese1995} and a refinement of \cite[Theorem $1.2$]{ABF} which, given an algebra whose congruence lattice has a pentagon as in \cref{fig:N5D1D2}, characterizes the sublattice of the congruence lattice of $\b$, considered as a subpower of the original algebra, generated using Freese's technique (\cref{thm:N5A(b)}). In the paper, we will produce explicit examples of algebras satisfying each of the cases in \cref{thm:N5A(b)}. As a consequence of \cref{thm:N5A(b)}, we also derived a characterization of congruence modular varieties in terms of lattice omissions (\cref{cor:N5char}). 

    In the last decade, Taylor varieties, i.e., varieties satisfying a non-trivial idempotent identity, have attracted significant attention for several reasons. Surprisingly, the class of Taylor varieties is a strong Mal'cev class \cite{Ols17}, and this particular class of varieties exhibits a profound connection with the Feder-Vardi conjecture, which was independently proven to be true in \cite{bulatov-dichotomy} and \cite{zhuk-dichotomy-short}. We apply Freese's technique to Taylor varieties (\cref{thm:K13D1}), proving that a variety is Taylor if and only if it omits $\alg D_{13}$, \cref{fig:D13}, a similar but different result with respect to \cite[Theorem 4.23]{KearnesKiss2013}. 

    We then focus on varieties satisfying a non-trivial congruence equation, which were classified in terms of lattice omission in \cite[Theorem $8.11$]{KearnesKiss2013}, and find a characterization of this class of algebras, \cref{thm:d2char}.

    The characterizations established in \cref{thm:N5A(b)}, \cref{thm:K13D1}, and \cref{thm:D2}, obtained using \emph{Freese's technique}, do more than provide descriptions of the aforementioned well-known classes of algebras—descriptions that can often also be obtained by other methods. Rather, they characterize when respectively a $\alg N_5$, $\alg D_1$, or $\alg D_2$ in a congruence lattice of a given algebra produces particular lattices in the congruence lattice of the algebra, whose underlying set is a selected congruence of the original algebra. The primary objective of this paper is, indeed, to further deepen our understanding of this technique, motivated by its proven effectiveness in the study of Mal'cev conditions.

    \subsection*{Organization of the paper} In Section~\ref{sec:notation} we introduce the necessary concepts and notation. Section~\ref{sec:FreeseTech} introduces the main methodology of the paper, \emph{Freese's technique}, with a simple example of its application and a list of results in the literature proven with this technique.  Next, we focus on the characterization of the admission of lattices using Freese's technique. In Section~\ref{sec:N5} we determine the lattices produced by a $\alg N_5$ as in \cref{fig:N5D1D2} via Freese's technique applied to $\b$ and then provide a characterization of congruence modular varieties, \cref{cor:N5char}. 
    In Section~\ref{sec:D1} we do the same for $\alg D_1$ (\cref{thm:K13D1}), providing a characterization of Taylor varieties \cref{TaylorThmorig}. Section~\ref{sec:D2} deals with $\alg D_2$, \cref{thm:D2}.  As an application, we present characterizations for the class of varieties satisfying a non-trivial congruence equation, \cref{thm:d2char}. The last Section~\ref{sec:Conclusions} is dedicated to ideas for future work and open problems.

	\section{Preliminaries and notation}\label{sec:notation}

    This section provides a review of some basic notions in lattice theory. For elementary concepts in general algebra, including lattices, algebras, and varieties, we refer to \cite{BurrisSanka}. For the general theory of Mal'cev conditions and Mal'cev classes,  we refer the reader to the classical treatment in \cite{Taylor1973} or the more contemporary approach presented in \cite{KearnesKiss2013}.

	For any variety $\vv V$, we introduce the following notation:
	$$
	\adm(\vv V) =\{\alg L \mid \alg L \leq \con(\alg A),  \alg A \in \vv  V\},
	$$
	and we say that $\vv V$ \emph{admits} $\alg L$ if $\alg L \in \adm(\vv V)$; otherwise, we say that $\vv V$ \emph{omits} $\alg L$.  The \emph{congruence variety} of $\vv V$ is the variety generated by $\adm(\vv V)$ and it is denoted by $\Con(\vv V)$; a variety of lattices is a congruence variety if it is equal to $\Con(\vv V)$ for some variety $\vv V$. In general, $\adm(\vv V) \subsetneq \Con(\vv V)$.
	
	Furthermore, we introduce the concept of \emph{filter of omission} or \emph{omission class}. Let $\Gamma$ be a class of lattices. Then
	$$
	\mathfrak F(\Gamma)=\{\vv V: \alg L \notin \adm(\vv V), \text{for all}\ \alg L \in \Gamma\}
	$$
	and we refer to $\mathfrak F(\Gamma)$ as the \emph{omission class of $\Gamma$}.
	We denote by $\NN_0 = \{0,1,\dots\} = \NN \cup\{0\}$ and $[n] = \{0,1,\dots, n\}$. Let $A$ be a set, and $S \subseteq A \times A$. Then we denote by $\Eq(S)$ the smallest equivalence relation over $A$ containing $S$. Furthermore, for an algebra $\alg A$ and a set of pairs $S \in A^2$, we denote by $\Cg(S)$ the congruence generated by $S$, i.e., the smallest congruence containing $S$.
	
	\section{Freese's  technique}\label{sec:FreeseTech}

	Let $\alg A$ be any algebra and $\a \in \con(\mathbf{A})$. Following the notation in \cite{ABF,Freese2024}, we denote by $\alg A^n(\a)$ the subalgebra of $\alg A^n$ whose universe is
	$$
	A^n(\a) :=\{(\vuc an) \mid (a_i,a_j) \in \a, i,j \le n\}.
	$$
    Let $\alg L$ be a given lattice, $\vv V$ be a variety of algebras, and $\alg A\in \vv V$. Our goal is to identify conditions on the congruence lattice of $\alg A$ that allow us to conclude that $\alg L \in \Con(\vv V)$. A stronger question is whether there exists $\a \in \con(\mathbf{A})$ such that $\alg L \le \op{Con}(\alg A^n(\a))$. The approach of studying sublattices of $\op{Con}(\alg A^n(\a))$ originates primarily from R. Freese, who has authored or coauthored nearly all published works on this topic. Following the notation in \cite{ABF}, we refer to this method as \emph{Freese's technique}. In this paper, we will consider the special case where $n=2$, omitting the superscript and writing $\alg A(\a) = \a$.
	
	The following is a list of results that can be, and have been, proven using Freese's technique.  Here, $\vv M_p$ denotes the congruence variety of vector spaces over a field of characteristic $p$.
	
	\begin{enumerate}
		\item There are non-modular varieties of lattices that are not congruence varieties (Nation \cite{Nation1974} without Freese's technique; then Freese, see \cite{Freese1995});
		\item every modular congruence variety $\vv V$ consists entirely of {\em arguesian} lattices, a stronger property with respect to modularity, (Freese and J\'onsson, \cite{FreeseJonsson1976}); hence the variety of all modular lattices is not a congruence variety,
		\item every non-modular congruence variety contains $\Con(\vv P)$, where $\vv P$ is Polin's variety (Freese and Day, \cite{DayFreese1980});
		\item every modular non-distributive congruence variety $\vv V$ contains $\vv M_p$, for some $p$ that is prime or $0$ (Freese, Herrmann and Huhn, \cite{FreeseHerrmanHuhn1981}).
	\end{enumerate}

    We also want to highlight a recent remarkable result obtained using Freese's technique for $n > 2$. This result was subsequently used in \cite[Theorem $1.2$]{Freese2024} to prove the following: If $\vv V$ is a variety such that $\Con(\vv V)$ satisfies a
non-trivial lattice identity, then $\Con(\vv V)$ is axiomatized by a finite set of equations
only if it is semidistributive.

    \begin{Thm}[Theorem 6.1 in \cite{Freese2024}]\label{freese}
        Let $\vv V$ be a variety having a weak difference term, and assume that $\vv V$ is not congruence meet semidistributive. Then there is a prime field $\alg P$ such that, for every finite $n$, the lattice of subspaces of a vector space over $\alg P$
of dimension $n$ lies in $\Con(\vv V)$.
    \end{Thm}
	
	For completeness, we include some fundamental results related to Freese's technique for $n =2$, which can also be found in \cite[Section $3$]{ABF}. For $n =2$, the technique was first introduced in \cite{FreeseJonsson1976} as \emph{duplication technique}. Note that for an algebra $\alg A$ with $\a \in \con(\alg A)$, we have that $\alg A(\a) \le \alg A^2$. In particular, $\alg A(\a)$ is a subdirect product of $\alg A^2$ and we can use the following notation for congruences on subdirect products. If $\th \in \con(\alg A)$, then:
	\begin{enumerate}
		\ib $\th_0 = \{( (a_0,a_1),(b_0,b_1)) \mid (a_0,b_0) \in \th\}$;
		\ib $\th_1 = \{( (a_0,a_1),(b_0,b_1)) \mid (a_1,b_1) \in \th\}$;
		\ib $\eta_0$,$\eta_1$ are the kernels of the projections of $\alg A(\a)$ on the two factors.
	\end{enumerate}
	By the Homomorphism  Theorem  \cite[Theorem $6.11$]{BurrisSanka} $\op{Con} (\alg A(\a)/\eta_0) \cong$  $ \op{Con} (\alg A(\a)/\eta_1)$ $ \cong \con (\alg A)$; moreover, the following facts are easy to prove  (see also \cite{FreeseJonsson1976}).

	\begin{Lem}\label{doublinglemma}  Let $\alg A$ be an algebra, $\g,\th \in \con(\alg A)$. Then the following holds in $\op{Con}(\alg A(\g))$:
		\begin{enumerate}
			\item if $\g \le \th$, then $\th_0=\th_1$;
			\item $\th_i = \eta_i \join (\th_0 \meet \th_1)$;
			\item $\g_0 = \eta_0 \join \eta_1$.
		\end{enumerate}
	\end{Lem}

    Throughout the paper, whenever we consider a sublattice of the congruence lattice of an algebra $\alg A$, we will assume that its bottom element coincides with the bottom element of $\con(\alg A)$. Although the results are stated in full generality and do not require this assumption, adopting it simplifies the exposition. Indeed, this entails no real loss of generality: if the bottom element of the chosen sublattice is a congruence $\theta \neq 0_{\alg A}$, one may pass to the quotient algebra $\alg A/\theta$, where the corresponding sublattice has a bottom element $0_{\alg A/\theta}$. After establishing the desired results in the quotient algebra, they can be transferred back to the original algebra via the Homomorphism Theorem. In particular, we simply note that, by the Homomorphism Theorem, given an algebra $\alg A$ and $\gamma \leq \beta$ in $\con(\alg A)$, then $(\alg A/\gamma)(\beta/\gamma) \cong \alg A(\beta)/(\gamma_0\wedge\gamma_1)$.
	
	\begin{figure}[hbtp]
		\begin{center}
			\begin{tikzpicture}[scale=0.8]
				\draw (0,0)-- (-1,1) -- (0,2) -- (0.6,1.4) -- (0.6,0.6) -- (0,0);
				\draw[fill] (0,0) circle [radius=0.05];
				\draw[fill] (-1,1) circle [radius=0.05];
				\draw[fill] (0,2) circle [radius=0.05];
				\draw[fill] (0.6,1.4) circle [radius=0.05];
				\draw[fill] (0.6,0.6) circle [radius=0.05];
				\node [right] at (0,-.1) {\footnotesize $0$};
				\node [right] at (0.6,0.6) {\footnotesize $\g$};
				\node [right] at (0.6,1.4) {\footnotesize $\a$};
				\node [left] at (-1,1) {\footnotesize $\beta$};
				\node [right] at (0,2.1) {\footnotesize $1$};
				
			\draw (10,2)-- (11,1) -- (10,0) -- (9,1)  -- (10,2) -- (10,1)--(10.5,0.5)-- (10,1)--(9.5,0.5); 
			\draw[fill] (10,1) circle [radius=0.05];
			\draw[fill] (10.5,0.5) circle [radius=0.05];
			\draw[fill] (9.5,0.5) circle [radius=0.05];
			\draw[fill] (10,0) circle [radius=0.05];
			\draw[fill] (9,1) circle [radius=0.05];
			\draw[fill] (11,1) circle [radius=0.05];
			\draw[fill] (10,2) circle [radius=0.05];
			\node [right] at (10,2.1) {\footnotesize $1$};
			\node [right] at (11,1){\footnotesize $\g$};
			\node [left] at (9,1){\footnotesize $\a$};
			\node [right] at (10,-0.1) {\footnotesize $0$};
			\node [left] at (9.5,0.4)  {\footnotesize $\mu$};
			\node [right] at (10.5,0.4) {\footnotesize $\delta$};
			\node [right] at (10,1.1)  {\footnotesize $\b$};

			\draw (5,2)-- (6,1) -- (5,0) -- (4,1)  -- (5,2) -- (4.5,1.5) -- (5, 1) -- (5.5,1.5) -- (5, 1) -- (5, 0); 
            \draw (10,2)-- (11,1) -- (10,0) -- (9,1)  -- (10,2) -- (10,1)--(10.5,0.5)-- (10,1)--(9.5,0.5); 
			\draw[fill] (5,1) circle [radius=0.05];
			\draw[fill] (5.5,1.5) circle [radius=0.05];
			\draw[fill] (4.5,1.5) circle [radius=0.05];
			\draw[fill] (5,0) circle [radius=0.05];
			\draw[fill] (4,1) circle [radius=0.05];
			\draw[fill] (6,1) circle [radius=0.05];
			\draw[fill] (5,2) circle [radius=0.05];
			\node [right] at (5,2.1) {\footnotesize $1$};
			\node [right] at (6,1){\footnotesize $\g$};
			\node [left] at (4,1){\footnotesize $\a$};
			\node [right] at (5,-0.1) {\footnotesize $0$};
			\node [left] at (4.5,1.6)  {\footnotesize $\mu$};
			\node [right] at (5.5,1.6) {\footnotesize $\delta$};
			\node [right] at (5,1)  {\footnotesize $\b$};
			\end{tikzpicture}
		\end{center}
		\caption{Respectively $\alg N_5$, $\alg D_1$ and $\alg D_2$}\label{fig:N5D1D2}
	\end{figure}

	As an example of basic application of Freese's technique, suppose that $\alg A$ is an algebra that is not congruence modular. Then $\alg N_5 \le \con(\alg A)$. If $\alg N_5$ is labeled as in Figure \ref{fig:N5D1D2}, let us apply Freese's technique to $\g$. Then, as a consequence of \cref{doublinglemma}, we can see that inside $\op{Con}(\alg A(\g))$ lies the sublattice $\alg L_{14}$ of Figure \ref{L15lat}.
	
	\begin{figure}[h]
		\begin{center}
			\begin{tikzpicture}
				\draw (0,0)-- (-1.5,1.5) -- (0,3) -- (1.5,1.5) -- (0,0) -- (0,1.5) -- (-1.5,3) -- (0,4.5) -- (1.5,3)-- (0,1.5);
				\draw (-1.5,1.5) -- (-1.5,3);
				\draw (0,3) -- (0,4.5);
				\draw (1.5,1.5) --(1.5,3);
				\draw[fill] (0,0) circle [radius=0.05];
				\draw[fill] (-1.5,1.5) circle [radius=0.05];
				\draw[fill] (0,3) circle [radius=0.05];
				\draw[fill] (1.5,1.5) circle [radius=0.05];
				\draw[fill] (0,1.5) circle [radius=0.05];
				\draw[fill] (-1.5,3) circle [radius=0.05];
				\draw[fill] (0,4.5) circle [radius=0.05];
				\draw[fill] (1.5,3) circle [radius=0.05];
				\draw[fill] (0,3.75) circle [radius=0.05];
				\node [right] at (0,-.1) {\footnotesize $0_{\alg A(\g)}$};
				\node [left] at (-1.5,1.5) {\footnotesize $\eta_0$};
				\node [right] at (1.5,1.5) {\footnotesize $\eta_1$};
				\node [left] at (-1.5,3) {\footnotesize $\beta_0$};
				\node [right] at (1.5,3) {\footnotesize $\beta_1$};
				\node [right] at (0,3) {\footnotesize $\g_0=\g_1$};
				\node [right] at (0,1.5) {\footnotesize $\beta_0 \meet \beta_1$};
				\node [right] at (0,3.75) {\footnotesize $\a_0=\a_1$};
				\node [right] at (0,4.5) {\footnotesize $1_0 = 1_1$};
			\end{tikzpicture}
		\end{center}
		\caption{$\alg L_{14}$}\label{L15lat}
	\end{figure}

	It follows that if $\vv V$ is not congruence modular then $\alg L_{14} \in \adm(\vv V)$ Moreover, since $\vv V(\alg N_5)$ is a proper subvariety of $\vv V(\alg L_{14})$, any variety $\vv U$
	of lattices such that $\alg N_5 \in \vv U$ but $\alg L_{14} \notin \vv U$ cannot be a congruence variety. This provides a proof that there exist non-modular varieties of lattices that are not congruence varieties \cite{FreeseJonsson1976}. Note that $\alg L_{14}$ is the sublattice generated by the two copies of $\alg N_5$ within $\con(\alg A(\g))$ and it is uniquely determined. R. Freese observed in \cite{Freese1995} that when the same duplication construction is applied to $\a$ or $\b$, then the sublattice generated by the two pentagons is no longer uniquely determined. An investigation of the second complex case can be found in \cite{ABF}.  

	\begin{Lem}[Lemma 3.2 of \cite{ABF}]\label{lemma3} Let $\alg A$ be an algebra such that $\con(\mathbf{A})$ has a sublattice as the $\alg N_5$ in Figure \ref{fig:N5D1D2}. Then the following hold in $\op{Con}(\alg A(\b))$:
		\begin{enumerate}
			\item $\g_0 \meet \g_1 \le \a_0 \meet \g_1, \g_0 \meet \a_1, \a_0 \meet \a_1$;
			\item $\a_0 \meet \a_1 \not\le \a_0 \meet \g_1, \g_0 \meet \a_1, \g_0 \meet \g_1$;
			\item  if $\a_0 \meet \g_1$ and $\g_0 \meet \a_1$ are comparable congruences, then
			$$
			\a_0 \meet \g_1 =\g_0 \meet \a_1= \g_0\meet\g_1;
			$$
			\item$\g_0 \join (\a_0 \meet \a_1) = \a_0$ and $\g_1 \join (\a_0 \meet \a_1) = \a_1.$
		\end{enumerate}
	\end{Lem}

	\begin{Thm}[Theorem 1.1 \cite{ABF}]\label{thm:N5old}
		Let $\alg A$ be an algebra with a non-modular congruence lattice. Then there is an $\alg N_5 \leq \con(\alg A)$ (labeled as in Figure \ref{fig:N5D1D2}) such that the sublattice $\alg L$ of $\con(\alg A(\beta))$ generated by $\{\alpha_0, \alpha_1, \gamma_0, \gamma_1, \beta_0\}$ is isomorphic to:
		\begin{enumerate}
			\item $\alg K$ (Figure \ref{M1K}) if and only if  $\alpha_0 \wedge \gamma_1$ and $\gamma_0 \wedge \alpha_1$ are not comparable;

            \item $\alg M_1$ (Figure \ref{M1K}) if and only if $\alpha_0 \wedge \gamma_1$ and $\gamma_0 \wedge \alpha_1$ are comparable.
		\end{enumerate}
	\end{Thm}

Note that for a non-congruence modular algebra $\alg A$, \cref{thm:N5old} tells that there exists a pentagon, labeled as in \cref{fig:N5D1D2}, within $\con(\alg A)$ such that if we consider an application of Freese's technique to $\b$, then the two copies of $\alg N_5$ in $\con(\alg A(\b))$ generate only two possible lattices.  However, this theorem does not provide any specific information regarding a fixed pentagon within $\con(\alg A)$. In the next section, we will show that the situation is much more complex in this case.

    \begin{figure}[htbp]
		\begin{center}
			\begin{tikzpicture}
            \draw (0,0) -- (-2,1.5) -- (-2,3) -- (0,4.5) -- (0,3) -- (0,4.5) -- (2,3) -- (2,1.5) -- (0,0);
				\draw (0,0) -- (0,1.5);
				\draw (-2,1.5) -- (0,3) -- (2,1.5) -- (0,0);
				\draw (-2,3) -- (0,1.5) -- (2,3);
				\draw (-2,2.2) -- (0,0.7) -- (2,2.2);
				\draw[fill] (0,0) circle [radius=0.05];
				\node [right] at (0,-.2) {\footnotesize $0_{\alg A(\b)}$};
				\draw[fill] (-2,1.5) circle [radius=0.05];
				\node [left] at (-2,1.5) {\footnotesize $\eta_0$};
				\draw[fill] (-2,3) circle [radius=0.05];
				\node [left] at (-2,3) {\footnotesize $\a_0$};
				\draw[fill] (0,4.5) circle [radius=0.05];
				\node [right] at (0,4.5) {\footnotesize $1_0 = 1_1$};
				\draw[fill] (0,3) circle [radius=0.05];
				\node [right] at (0,3) {\footnotesize $\b_0=\b_1$};
				\draw[fill] (2,3) circle [radius=0.05];
				\node [right] at (2,3) {\footnotesize $\a_1$};
				\draw[fill] (2,1.5) circle [radius=0.05];
				\node [right] at (2,1.4) {\footnotesize $\eta_1$};
				\draw[fill] (0,1.5) circle [radius=0.05];
				\node [right] at (0,1.5) {\footnotesize $\a_0\meet \a_1$};
				\draw[fill] (-2,2.2) circle [radius=0.05];
				\node [left] at (-2,2.2) {\footnotesize $\g_0$};
				\draw[fill] (0,0.7) circle [radius=0.05];
				\node [right] at (0,0.7) {\footnotesize $\g_0\meet \g_1$};
				\draw[fill] (2,2.2) circle [radius=0.05];
				\node [right] at (2,2.2) {\footnotesize $\g_1$};

            
				\draw (6,0) -- (4,1.5) -- (4,3.5) -- (6,5) -- (6,3) -- (6,5.) -- (8,3.5) -- (8,1.5) -- (6,0);
				\draw (6,0) -- (6,0.7);
				\draw (4,1.5) -- (6,3) -- (8,1.5) -- (6,0);
				\draw (4,3.5) -- (6,2)-- (6,1.5) -- (6,2) -- (8,3.5);
				\draw (4,2.2) -- (6,0.7) -- (8,2.2);
				\draw (6,1.5) -- (6.48,1.05);
				\draw (6,1.5) -- (5.52,1.05);
				\draw[fill] (6,0) circle [radius=0.05];
				\node [right] at (6,-.2) {\footnotesize $0_{\alg A(\b)}$};
				\draw[fill] (4,1.5) circle [radius=0.05];
				\node [left] at (4,1.5) {\footnotesize $\eta_0$};
				\draw[fill] (4,3.5) circle [radius=0.05];
				\node [left] at (4,3.5) {\footnotesize $\a_0$};
				\draw[fill] (6,5) circle [radius=0.05];
				\node [right] at (6,5) {\footnotesize $1_0 = 1_1$};
				\draw[fill] (6,3) circle [radius=0.05];
				\node [right] at (6,3) {\footnotesize $\b_0=\b_1$};
				\draw[fill] (8,3.5) circle [radius=0.05];
				\node [right] at (8,3.5) {\footnotesize $\a_1$};
				\draw[fill] (8,1.5) circle [radius=0.05];
				\node [right] at (8,1.4) {\footnotesize $\eta_1$};
				\draw[fill] (6,1.5) circle [radius=0.05];
				\node [right] at (6,1.5) {\footnotesize $\th$ };
				\draw[fill] (4,2.2) circle [radius=0.05];
				\node [left] at (4,2.2) {\footnotesize $\g_0$};
				\draw[fill] (6,0.7) circle [radius=0.05];
				\node [right] at (6,0.7) {\footnotesize $\g_0\meet \g_1$};
				\draw[fill] (8,2.2) circle [radius=0.05];
				\node [right] at (8,2.2) {\footnotesize $\g_1$};
				\draw[fill] (6,2) circle [radius=0.05];
				\node [right] at (6,2) {\footnotesize $\a_0\meet \a_1$};
				\draw[fill] (5.52,1.05) circle [radius=0.05];
				\node [left] at (5.52,1.05) {\footnotesize $\g_0 \meet \a_1$};
				\draw[fill] (6.48,1.05) circle [radius=0.05];
				\node [right] at (6.48,1.05) {\footnotesize $\a_0\meet \g_1$};
			\end{tikzpicture}
		\end{center}
		\caption{Respectively, the lattices $\alg M_1$ and $\alg K$}\label{M1K}
	\end{figure}

	\section{An application of Freese's technique to $\alg N_5$}\label{sec:N5}
	
	In this section, our aim is to apply Freese's technique to the lattice $\alg N_5$ to investigate the shape of the congruence lattice $\con(\b)$. This lattice was previously studied in \cite{ABF}, where a partial characterization of the sublattice generated by $\{\alpha_0, \alpha_1, \gamma_0, \gamma_1, \beta_0\}$ in $\con(\b)$ was provided. That description was obtained under the assumption that $[\g,\a] =\{\g,\a\}$. Here, we further refine this result by removing this hypothesis and introducing a family of lattices whose intersection of omission filters is equivalent to the filter of omission of $\alg N_5$. 
    
    Let $\alg N_5$ be the lattice depicted in \cref{fig:N5D1D2}. In \cref{fig:KiMi} and \cref{fig:Kinf}  we introduce two families of lattices, $\{\alg K_i\}_{i \in \mathbb{N}}$ and $\{\alg M_i\}_{i \in \mathbb{N}}$, along with the lattice $\alg K_{\infty}$ (\cref{fig:Kinf}), where $\theta_i = (\g_0^i \wedge \g_1^{i+1}) \vee (\g_0^{i+1} \wedge \g_1^{i})$, for all $i \in I$. These lattices are produced using Freese's technique applied to $\b$ of $\alg N_5$. Note that we use the superscript index instead of the usual underscript since  the latter would be in contrast with the usual notation of non-skew congruences.

 \begin{figure}[htbp]
		\begin{center}
			\begin{tikzpicture}[scale=0.70]
				
				\draw (-3,2) --(0,0) -- (3,2);
				\draw (-3,3) --(0,1) -- (3,3);
				\draw (-3,5) --(0,3) -- (3,5);

				\draw (-3,7) --(0,5) -- (3,7);
				\draw (-3,10) --(0,8) -- (3,10);
				\draw (3,2) --(0,10)--(-3,2);
				\draw (3,10)--(0,12)--(-3,10);
				\draw (3,8)--(0,6)--(-3,8);
				
				\draw (0,0) --(0,1);
				\draw (0,2) --(0,3);
				\draw (3,2) --(3,5);
				\draw (-3,2) --(-3,5);
				\draw (3,5) --(3,7);
				\draw (-3,5) --(-3,7);
				\draw (0,4) -- (0,5);
				\draw (0,10) --(0,12);
                \draw (0,7) --(0,8);
                \draw (-3,8) --(-3,10);
                \draw (3,8) --(3,10);

				\draw [dashed](-3,7) --(-3,8);
				\draw [dashed](3,7) --(3,8);
				\draw [dashed](0,5) --(0,6);
				
				\draw(-0.75,1.5) -- (0,2) -- (0.75,1.5);
				\draw(-0.75,3.5) -- (0,4) -- (0.75,3.5);
				\draw(-0.75,6.5) -- (0,7) -- (0.75,6.5);

				\draw[fill] (0,0) circle [radius=0.05];
				\node [left] at (0,0) {\footnotesize $0_{\alg A(\b)}$};
				\draw[fill] (0,1) circle [radius=0.05];
				\node [right] at (0,1) {\tiny $\g^0_0 \wedge \g_0^1$};
				\draw[fill] (0,2) circle [radius=0.05];
				\node [right] at (0,2) {\tiny $\theta_0$};
				\draw[fill] (0,3) circle [radius=0.05];
				\node [right] at (0,3) {\tiny $\g^1_0 \wedge \g_1^1$};
				\draw[fill] (0,4) circle [radius=0.05];
				\node [right] at (0,4) {\tiny $\theta_1$};
                \draw[fill] (0,7) circle [radius=0.05];
                \node [right] at (0,7) {\tiny $\theta_i$};
				\draw[fill] (0,5) circle [radius=0.05];
				\node [right] at (0,5) {\tiny $\g^{2}_0 \wedge \g_1^{2}$};
				\node [right] at (0,6) {\tiny $\g^i_0 \wedge \g^i_1$};
				\draw[fill] (0,6) circle [radius=0.05];
                \node [right] at (0,8) {\tiny $\a_0 \wedge \a_1$};
				\draw[fill] (0,8) circle [radius=0.05];
				\node [right] at (0,10) {\tiny $\b_0 = \b_1$};
				\draw[fill] (0,10) circle [radius=0.05];
				\node [right] at (0,12) {\tiny $1_0 = 1_1$};
				\draw[fill] (0,12) circle [radius=0.05];

				\draw[fill] (-0.75,1.5) circle [radius=0.05];
				\node [left] at (-0.85,1.5) {\tiny $\g^{0}_0 \wedge \g_1^{1}$};
				\draw[fill] (0.75,1.5) circle [radius=0.05];
				\node [right] at (0.85,1.5) {\tiny $\g^{1}_0 \wedge \g_1^{0}$};
				\draw[fill] (-0.75,3.5) circle [radius=0.05];
				\node [left] at (-0.85,3.5) {\tiny $\g^{1}_0 \wedge \g_1^{2}$};
				\draw[fill] (0.75,3.5) circle [radius=0.05];
				\node [right] at (0.85,3.5) {\tiny $\g^{2}_0 \wedge \g_1^{1}$};

                \draw[fill] (-0.75,6.5) circle [radius=0.05];
				\node [left] at (-0.85,6.5) {\tiny $\g^{i}_0 \wedge \a_1$};
				\draw[fill] (0.75,6.5) circle [radius=0.05];
				\node [right] at (0.85,6.5) {\tiny $\a_0 \wedge \g_1^{i}$};

				\draw[fill] (3,2) circle [radius=0.05];
				\node [right] at (3,2) {\tiny $\eta_1$};
				\draw[fill] (3,3) circle [radius=0.05];
				\node [right] at (3,3) {\tiny $\g^0_1$};
				\draw[fill] (3,5) circle [radius=0.05];
				\node [right] at (3,5) {\tiny $\g^1_1$};
				\draw[fill] (3,7) circle [radius=0.05];
				\node [right] at (3,7) {\tiny $\g^{2}_1$};
				\draw[fill] (3,8) circle [radius=0.05];
				\node [right] at (3,8) {\tiny $\g^i_1$};
                \draw[fill] (3,10) circle [radius=0.05];
				\node [right] at (3,10) {\tiny $\a_1$};

				\draw[fill] (-3,2) circle [radius=0.05];
				\node [left] at (-3,2) {\tiny $\eta_0$};
				\draw[fill] (-3,3) circle [radius=0.05];
				\node [left] at (-3,3) {\tiny $\g^0_0$};
				\draw[fill] (-3,5) circle [radius=0.05];
				\node [left] at (-3,5) {\tiny $\g^1_0$};
				\draw[fill] (-3,7) circle [radius=0.05];
				\node [left] at (-3,7) {\tiny $\g^{2}_0$};
				\draw[fill] (-3,8) circle [radius=0.05];
				\node [left] at (-3,8) {\tiny $\g^i_0$};
                \draw[fill] (-3,10) circle [radius=0.05];
				\node [left] at (-3,10) {\tiny $\a_0$};


                \draw (5,2) --(8,0) -- (11,2);
				\draw (5,3) --(8,1) -- (11,3);
				\draw (5,5) --(8,3) -- (11,5);

				\draw (5,7) --(8,5) -- (11,7);
				\draw (5,10) --(8,8) -- (11,10);
				\draw (5,2) --(8,10)--(11,2);
				\draw (11,10)--(8,12)--(5,10);
				\draw (11,8)--(8,6)--(5,8);
				
				\draw (8,0) --(8,1);
				\draw (8,2) --(8,3);
				\draw (11,2) --(11,5);
				\draw (5,2) --(5,5);
				\draw (11,5) --(11,7);
				\draw (5,5) --(5,7);
				\draw (8,4) -- (8,5);
				\draw (8,10) --(8,12);
                \draw (8,6) --(8,8);
                \draw (5,8) --(5,10);
                \draw (11,8) --(11,10);

				\draw [dashed](5,7) --(5,8);
				\draw [dashed](11,7) --(11,8);
				\draw [dashed](8,5) --(8,6);
				
				\draw(7.25,1.5) -- (8,2) -- (8.75,1.5);
				\draw(7.25,3.5) -- (8,4) -- (8.75,3.5);

				\draw[fill] (8,0) circle [radius=0.05];
				\node [left] at (8,0) {\footnotesize $0_{\alg A(\b)}$};
				\draw[fill] (8,1) circle [radius=0.05];
				\node [right] at (8,1) {\tiny $\g^0_0 \wedge \g_0^1$};
				\draw[fill] (8,2) circle [radius=0.05];
				\node [right] at (8,2) {\tiny $\theta_0$};
				\draw[fill] (8,3) circle [radius=0.05];
				\node [right] at (8,3) {\tiny $\g^1_0 \wedge \g_1^1$};
				\draw[fill] (8,4) circle [radius=0.05];
				\node [right] at (8,4) {\tiny $\theta_1$};
				\draw[fill] (8,5) circle [radius=0.05];
				\node [right] at (8,5) {\tiny $\g^{2}_0 \wedge \g_1^{2}$};
				\node [right] at (8,6) {\tiny $\g^i_0 \wedge \g^i_1$};
				\draw[fill] (8,6) circle [radius=0.05];
                \node [right] at (8,8) {\tiny $\a_0 \wedge \a_1$};
				\draw[fill] (8,8) circle [radius=0.05];
				\node [right] at (8,10) {\tiny $\b_0 = \b_1$};
				\draw[fill] (8,10) circle [radius=0.05];
				\node [right] at (8,12) {\tiny $1_0 = 1_1$};
				\draw[fill] (8,12) circle [radius=0.05];

				\draw[fill] (7.25,1.5) circle [radius=0.05];
				\node [left] at (7.15,1.5) {\tiny $\g^{0}_0 \wedge \g_1^{1}$};
				\draw[fill] (8.75,1.5) circle [radius=0.05];
				\node [right] at (8.85,1.5) {\tiny $\g^{1}_0 \wedge \g_1^{0}$};
				\draw[fill] (7.25,3.5) circle [radius=0.05];
				\node [left] at (7.15,3.5) {\tiny $\g^{1}_0 \wedge \g_1^{2}$};
				\draw[fill] (8.75,3.5) circle [radius=0.05];
				\node [right] at (8.85,3.5) {\tiny $\g^{2}_0 \wedge \g_1^{1}$};

				\draw[fill] (11,2) circle [radius=0.05];
				\node [right] at (11,2) {\tiny $\eta_1$};
				\draw[fill] (11,3) circle [radius=0.05];
				\node [right] at (11,3) {\tiny $\g^0_1$};
				\draw[fill] (11,5) circle [radius=0.05];
				\node [right] at (11,5) {\tiny $\g^1_1$};
				\draw[fill] (11,7) circle [radius=0.05];
				\node [right] at (11,7) {\tiny $\g^{2}_1$};
				\draw[fill] (11,8) circle [radius=0.05];
				\node [right] at (11,8) {\tiny $\g^i_1$};
                \draw[fill] (11,10) circle [radius=0.05];
				\node [right] at (11,10) {\tiny $\a_1$};

				\draw[fill] (5,2) circle [radius=0.05];
				\node [left] at (5,2) {\tiny $\eta_0$};
				\draw[fill] (5,3) circle [radius=0.05];
				\node [left] at (5,3) {\tiny $\g^0_0$};
				\draw[fill] (5,5) circle [radius=0.05];
				\node [left] at (5,5) {\tiny $\g^1_0$};
				\draw[fill] (5,7) circle [radius=0.05];
				\node [left] at (5,7) {\tiny $\g^{2}_0$};
				\draw[fill] (5,8) circle [radius=0.05];
				\node [left] at (5,8) {\tiny $\g^i_0$};
                \draw[fill] (5,10) circle [radius=0.05];
				\node [left] at (5,10) {\tiny $\a_0$};

			\end{tikzpicture}
		\end{center}
		\caption{$\alg K_{i+1}$ and $\alg M_{i+1}$}\label{fig:KiMi}
	\end{figure}

    \begin{figure}[htbp]
		\begin{center}
			\begin{tikzpicture}[scale=0.90]
				
				\draw (-3,2) --(0,0) -- (3,2);
				\draw (-3,3) --(0,1) -- (3,3);
				\draw (-3,5) --(0,3) -- (3,5);

				\draw (-3,7) --(0,5) -- (3,7);
				\draw (-3,9) --(0,7) -- (3,9);
				\draw (3,2) --(0,9)--(-3,2);
				\draw (3,9)--(0,11)--(-3,9);
				
				\draw (0,0) --(0,1);
				\draw (0,2) --(0,3);
				\draw (3,2) --(3,5);
				\draw (-3,2) --(-3,5);
				\draw (3,5) --(3,7);
				\draw (-3,5) --(-3,7);
				\draw (0,4) -- (0,5);
				\draw (0,9) --(0,11);

				\draw [dashed](-3,7) --(-3,9);
				\draw [dashed](3,7) --(3,9);
				\draw [dashed](0,5) --(0,7);

				\draw(-0.75,1.5) -- (0,2) -- (0.75,1.5);
				\draw(-0.75,3.5) -- (0,4) -- (0.75,3.5);

				\draw[fill] (0,0) circle [radius=0.05];
				\node [left] at (0,0) {\footnotesize $0_{\alg A(\b)}$};
				\draw[fill] (0,1) circle [radius=0.05];
				\node [right] at (0,1) {\tiny $\g^0_0 \wedge \g_0^1$};
				\draw[fill] (0,2) circle [radius=0.05];
				\node [right] at (0,2) {\tiny $\theta_0$};
				\draw[fill] (0,3) circle [radius=0.05];
				\node [right] at (0,3) {\tiny $\g^1_0 \wedge \g_1^1$};
				\draw[fill] (0,4) circle [radius=0.05];
				\node [right] at (0,4) {\tiny $\theta_1$};
				\draw[fill] (0,5) circle [radius=0.05];
				\node [right] at (0,5) {\tiny $\g^{2}_0 \wedge \g_1^{2}$};
                \node [right] at (0,7) {\tiny $\a_0 \wedge \a_1$};
				\draw[fill] (0,7) circle [radius=0.05];
				\node [right] at (0,9) {\tiny $\b_0 = \b_1$};
				\draw[fill] (0,9) circle [radius=0.05];
				\node [right] at (0,11) {\tiny $1_0 = 1_1$};
				\draw[fill] (0,11) circle [radius=0.05];

				\draw[fill] (-0.75,1.5) circle [radius=0.05];
				\node [left] at (-0.85,1.5) {\tiny $\g^{0}_0 \wedge \g_1^{1}$};
				\draw[fill] (0.75,1.5) circle [radius=0.05];
				\node [right] at (0.85,1.5) {\tiny $\g^{1}_0 \wedge \g_1^{0}$};
				\draw[fill] (-0.75,3.5) circle [radius=0.05];
				\node [left] at (-0.85,3.5) {\tiny $\g^{1}_0 \wedge \g_1^{2}$};
				\draw[fill] (0.75,3.5) circle [radius=0.05];
				\node [right] at (0.85,3.5) {\tiny $\g^{2}_0 \wedge \g_1^{1}$};

				\draw[fill] (3,2) circle [radius=0.05];
				\node [right] at (3,2) {\tiny $\eta_1$};
				\draw[fill] (3,3) circle [radius=0.05];
				\node [right] at (3,3) {\tiny $\g^0_1$};
				\draw[fill] (3,5) circle [radius=0.05];
				\node [right] at (3,5) {\tiny $\g^1_1$};
				\draw[fill] (3,7) circle [radius=0.05];
				\node [right] at (3,7) {\tiny $\g^{2}_1$};
                \draw[fill] (3,9) circle [radius=0.05];
				\node [left] at (3,9) {\tiny $\a_1$};

				\draw[fill] (-3,2) circle [radius=0.05];
				\node [left] at (-3,2) {\tiny $\eta_0$};
				\draw[fill] (-3,3) circle [radius=0.05];
				\node [left] at (-3,3) {\tiny $\g^0_0$};
				\draw[fill] (-3,5) circle [radius=0.05];
				\node [left] at (-3,5) {\tiny $\g^1_0$};
				\draw[fill] (-3,7) circle [radius=0.05];
				\node [left] at (-3,7) {\tiny $\g^{2}_0$};
                \draw[fill] (-3,9) circle [radius=0.05];
				\node [left] at (-3,9) {\tiny $\a_0$};
    \end{tikzpicture}
		\end{center}
		\caption{$\alg K_{\infty}$}\label{fig:Kinf}
	\end{figure}

We will prove that the omission of these two families of lattices and of $\alg K_{\infty}$ is equivalent to the omission of $\alg N_5$ using Freese's technique and that the presence of an occurrence of these lattices in $\con(\alg A(\b))$ is strictly connected with the behavior of the relational product of $\b$ and $\gamma$. Specifically, from \cite[Theorem $2$]{ABF}, we know that if an algebra $\alg A$ has an $\alg N_5$ as in Figure \ref{fig:N5D1D2} as a sublattice of $\con(\alg A)$  and $[\gamma, \alpha]$ = $\{\gamma, \alpha\}$, then the sublattice $\alg L$ of $\con(\alg A(\beta))$ generated by $\{\alpha_0, \alpha_1, \gamma_0, \gamma_1, \beta_0\}$ is either isomorphic to $\alg K$ if  $\beta \circ \gamma \circ \beta \cap \alpha \not\subseteq \gamma $, or to $\alg M_1$ otherwise.

We will improve this result by removing the assumption $[\gamma, \alpha]$ = $\{\gamma, \alpha\}$, thereby providing a general characterization of the lattice generated by $\{\alpha_0, \alpha_1,$ $ \gamma_0, \gamma_1, \beta_0\}$. To achieve this, we first show several auxiliary results that we need for the proof of the main theorem in this section.  First of all, we generalize \cite[Lemma 3.4]{ABF} observing that its proof remains valid under weaker assumptions.

\begin{Lem}[Lemma 3.4 of \cite{ABF}]\label{lemma1} Let $\alg A$ be an algebra such that $\a,\b,\g \in \con(\alg A)$, such that $\a \wedge \b = 0_{\alg A}$, and $\g < \a$. Then in $\con(\alg A(\beta))$ we have:
		$$
		(\a_0 \meet \g_1) \join (\g_0 \meet \a_1) < \a_0 \meet \a_1.
		$$
	\end{Lem}

    \begin{proof} Let $\alpha, \gamma$ be as in the hypothesis. We will show that for any $(a,b) \in \a\setminus \g$,
		$$
		((a,a),(b,b)) \notin (\a_0 \meet \g_1) \join (\g_0 \meet \a_1);
		$$
		this is enough to prove the thesis.
		
		Suppose that $((a,a), (b,b)) \in (\a_0 \meet \g_1) \join (\g_0 \meet \a_1)$; then there is an $n \in \mathbb N$ and $(u_0,v_0), \dots,(u_n,v_n) \in \b$ such that $(a,a) = (u_0,v_0)$, $(b,b)  = (u_n, v_n)$ and
		\begin{align*}
			&(u_i,v_i) \mathrel{\a_0 \meet \g_1} (u_{i+1},v_{i+1}) \qquad\text{for $i$ even}\\
			&(u_i,v_i) \mathrel{\g_0 \meet \a_1} ( u_{i+1},v_{i+1}) \qquad\text{for $i$ odd}.
		\end{align*}
		Thus we have the following relations
		\begin{center}
			\begin{tikzcd}
				a \arrow[dd,dash]\arrow[r,"\a", dash]&u_1\arrow[dd,"\b",dash]\arrow[r,"\g", dash]&u_2\arrow[dd,"\b",dash]\arrow[r,"\a", dash]&u_3\arrow[dd,"\b",dash]\arrow[r,dash,dotted]&\arrow[r,dash,dotted]&u_{n-1}\arrow[dd,"\b",dash]\arrow[r,"\a", dash]&b\arrow[dd,dash] \\
				\\
				a\arrow [r,"\g",dash]&v_1\arrow[r,"\a", dash]&v_2\arrow[r,"\g", dash]&v_3\arrow[r,dash,dotted]&\arrow[r,dash,dotted]&v_{n-1}\arrow[r,"\g", dash]&b
			\end{tikzcd}
		\end{center}
		Since $\g \le \a$, we get $(u_1,v_1) \in \a$ and thus $(u_1,v_1)  \in \a \meet \b = 0_\alg A$, hence $u_1=v_1$. Continuing this argument, we get $u_i=v_i$ for $i \le n$, so
		$$
		a \mathrel{\g} v_1=u_1 \mathrel{\g} u_2=v_2 \mathrel{\g} \dots \mathrel{\g} b.
		$$
		Hence $(a,b) \in \g$ that is a contradiction. From this, we can conclude that $((a,a),( b,b)) \notin (\a_0 \meet \g_1) \join (\g_0 \meet \a_1)$ and thus $(\a_0 \meet \g_1) \join (\g_0 \meet \a_1) < \a_0 \meet \a_1$.
	\end{proof}

    We included the proof of the previous lemma for completeness since the hypothesis is slightly different with respect to \cite[Lemma 3.4]{ABF}.

 \begin{Lem}\label{lem:skew}
     Let $\alg A$ be an algebra and let $\a,\g,\b \in \con(\alg A)$. Then the following are equivalent:
     \begin{enumerate}
         \item $(\b\ \circ\ \g\ \circ \b) \cap \a \not\subseteq \g$;
         \item $\g_0 \wedge \g_1 < \a_0 \wedge \g_1$ in $\con(\b)$.
     \end{enumerate}
 \end{Lem}

 \begin{proof}
	
        The proof is essentially the same of \cite[Theorem 1.2]{ABF}, but for sake of completeness, we include it. For $(1) \implies(2)$, let $(a,b) \in (\beta \circ^3 \gamma) \cap \alpha \setminus \gamma$. Then there exist $x_1, x_2 \in A$ such that $a \ \beta\ x_1\ \gamma\ x_2\ \beta\ b$. Thus, we have that the following relation holds:
		\begin{center}
			\begin{tikzcd}
				a \arrow[dd,"\b",dash]\arrow[r,"\a", dash]&b\arrow[dd,"\b",dash]
				\\
				\\
				x_1\arrow [r,"\g",dash]&x_2
			\end{tikzcd}
		\end{center}
		
		Hence, $((a,x_1), (b, x_2)) \in \alpha_0 \wedge \gamma_1$ and thus $\alpha_0 \wedge \gamma_1 > \gamma_0 \wedge \gamma_1$. 
	
        $(2) \implies(1)$ follows directly from the definition of $\a_0 \wedge \g_1$ and from the hypothesis.
	\end{proof}

 \begin{Lem}\label{lem:etameet}
     Let $\alg A$ be an algebra and let $\a,\b \in \con(\alg A)$ such that $\a \wedge \b = 0_{\alg A}$. Then $\a_i \wedge \eta_j = 0_{\alg A(\b)}$ for $\{i,j\} = \{0,1\}$.
 \end{Lem}

 \begin{proof}
     Let us fix $i = 0$ and $j = 1$. Then $((a,x_1), (b, x_2)) \in \a_0 \wedge \eta_1$ if and only if $(a,b) \in \a$ and $x_1 = x_2 = x$. Thus, we are in the following situation:
     \begin{center}
			\begin{tikzcd}
				a \arrow[dd,"\b",dash]\arrow[r,"\a", dash]&b\arrow[dd,"\b",dash]
				\\
				\\
				x\arrow [r,"0_{\alg A}",dash]&x
			\end{tikzcd}
		\end{center}
    which implies $a = b$ since $\a \wedge \b = 0_{\alg A}$.
 \end{proof}

In order to generalize \cite[Theorem $2$]{ABF} we introduce a family of congruences that allows to further refine this result.

\begin{Def}\label{def:cong_chain}
	Let $\alg A$ be an algebra and let $\con(\alg A)$ be its congruence lattice. Let $\a, \b, \g \in \con(\alg A)$ be such that $\g \leq \a$. Then we define the family of congruences $\{\g^i(\alpha,\b,\gamma)\}_{i \in \NN_0}$ such that $\gamma^0(\alpha,\b,\gamma) = \gamma$, $\gamma^{i + 1}(\alpha,\b,\gamma) =  \Cg((\b \circ \g^i \circ \b)\cap \a)$. 
\end{Def}

Henceforth, we will omit the dependency of $\gamma^i(\alpha,\b,\gamma)$ writing only $\gamma^i$ when the three congruences are clear in the context. Moreover, we observe that $ \g^i \leq \gamma^{i+1} \leq \a$, for all $i \in \NN_0$. We note that $(\b \circ \g^i \circ \b) \cap \a$ is a tolerance (reflexive and symmetric subalgebra of $\alg A^2$) and therefore the congruence generated by the tolerance $(\b \circ \g^i \circ \b) \cap \a$ is its transitive closure. The family of congruences $\{\g_i\}_{i \in \NN_0}$ forms a (potentially infinite) ascending chain lying in the interval $[\g,\a]$. Furthermore, if $\b \circ \g^i \circ \b \cap \a = \g^i$, then the ascending chain stabilizes and $\gamma^i = \g^{i+1}$.

  \begin{Lem}\label{lem:gen_gi}
    Let $\alg A$ be an algebra and let $\a,\g,\b \in \con(\alg A)$ with $\g \leq \a$ and let $i \in \NN_0$. Then:
    \begin{enumerate}
        \item $\g^{i+1}_0$ is in the sublattice of $\con(\alg A (\b))$ generated by $\{ \eta_0, \a_0, \g^i_1\}$;
        \item $\g^i <\g^{i+1}$ if and only if $\g^i  \subset (\b\ \circ\ \g^i\ \circ \b) \cap \a$. 
    \end{enumerate}
 \end{Lem}

 \begin{proof}
     
    $(1)$ basically follows from \cref{def:cong_chain}, since $((a,x_1),(b,x_2)) \in \g^{i+1}_0$ if and only if there exist $u_0,v_0\dots,u_n,v_{2n} \in \alg A$ such that $(u_0,v_0) = (a,x_1)$, $(u_n,v_{2n}) = (b,x_2)$ and

     \begin{center}
			\begin{tikzcd}
				a \arrow[dd,dash]\arrow[r, dash]&u_0\arrow[dd,"\b",dash]\arrow[r,"\a", dash]&u_1\arrow[dd,"\b",dash]\arrow[r, dash]&u_1\arrow[dd,"\b",dash]\arrow[r,dash,dotted]&\arrow[r,dash,dotted]&u_{n}\arrow[dd,"\b",dash]\arrow[r, dash]&b\arrow[dd,dash] \\
				\\
				x_1&v_1\arrow[r,"\g^{i}", dash]&v_2 &v_3\arrow[r,dash,dotted]&\arrow[r,dash,dotted]&v_{2n-1}&x_2.
			\end{tikzcd}
		\end{center}

        Thus, $\g_0^{i+1} = \eta_0 \vee (\a_0 \wedge \g_1^i)$. Moreover, $(2)$ follows from \cref{def:cong_chain}.
 \end{proof}

We are now ready to prove the main statement of the section.

 \begin{Thm}\label{thm:N5A(b)}
    Let $\alg A$ be an algebra and let $\alg N_5 \leq \con(\alg A)$ labeled as in \cref{fig:N5D1D2}. Then, for all $i \in \mathbb{N}_0$, the sublattice of $\con(\alg A(\b))$ generated by $\{\a_0, \a_1,  \b_0, \g_0^0,$ $ \g_1^0\}$ is isomorphic to:
    \begin{enumerate}
        \item $\alg M_{i+1}$ if and only if $\g^{i-1} <\g^{i+1} = \g^i <\a$;
        \item $\alg K_{i+1}$ if and only if $\g^{i} < \g^{i+1} =\a$;
        \item $\alg K_{\infty}$ if and only if $\g^{i} < \g^{i+1}$, for all $i \in \mathbb{N}_0$.
    \end{enumerate}
 \end{Thm}

 \begin{proof}
     First, as observed in \cref{sec:FreeseTech}, we may assume without loss of generality that the bottom element of $\alg N_5$ coincides with the bottom element of $\con(\alg A)$. Then, we make some observations about the common structure shared by all the cases. We note that $\op{Con}(\alg A(\b)/\eta_s)$ is isomorphic to $\con(\alg A)$ through the natural isomorphism $f(\lambda) = \lambda_s$, for $s \in \{0,1\}$. Furthermore, if $\g^j < \g^{j+1}$, then also $\g^j_0 \wedge \g^j_1  < \g^{j+1}_0 \wedge \g^{j+1}_0$, for $j \in \mathbb N_0$. Moreover, by \cref{lem:skew} and \cref{lem:gen_gi}, whenever $\g^j <\g^{j+1}$, we have that 
     \[
     \g^j_0 \wedge \g^{j}_1 <\g^j_0 \wedge \g^{j+1}_1, \g^{j+1}_0 \wedge \g^{j}_1 < \g^{j+1}_0 \wedge \g^{j+1}_1.
     \]

   Hence, by \cref{lem:gen_gi}, we can observe that

     \[
     \eta_0 \vee (\a_0 \wedge \g_1^j) = \eta_0 \vee \g^j_0 \vee (\a_0 \wedge \g_1^j) = \g^{j+1}_0,
     \] 

     \noindent for all $j \in \mathbb{N}_0$. Moreover, $\g_0^{j+1} \wedge \g^j_1 = \a_0 \wedge \g^j_1$, for all $j \in \mathbb{N}_0$, follows from the \cref{def:cong_chain} since $((a,b),(x_1,x_2)) \in \a_0 \wedge \g^j_1$ implies $(a,b) \in \g^{j+1}$ and $(x_1,x_2) \in \g^j$. The other joins and meets follow from the order relations between the factor congruences.

     What remains to prove is basically which of the congruences in $\alg K_{\infty}$ is collapsing. By \cref{lemma1}, $\theta_j = (\g^{j+1}_0 \meet \g^j_1) \join (\g_0^j \meet \g^{j+1}_1) < \g^{j+1}_0 \meet \g^{j+1}_1$, for all $j \in \mathbb N_0$, with $\g^{i} <\g^{i+1}$.

     Case $\g^{i-1} <\g^{i+1} = \g^i <\a$: then, by \cref{lem:skew}, $\g^i_0 \wedge \a_1 = \g^i_0 \wedge \g^i_1$, the chain $\{\g^j\}_{j \in \mathbb N_0}$ collapses for $j \geq i$, and thus the sublattice of $\con(\alg A(\b))$ generated by $\{\a_0, \a_1,  \b_0, \g_0^0, \g_1^0\}$ is isomorphic to $\alg M_{i+1}$.

     Case $\g^{i} <\g^{i+1} =\a$: then, by \cref{lem:skew}, $\g^i_0 \wedge \a > \g^i_0 \wedge \g^i_1$, $\g_0^i \vee(\g^i_0 \wedge \a_1) = \a_0$, and the chain $\{\g^j\}_{j \in \mathbb N_0}$ collapses for $j > i$. Furthermore, by \cref{lemma1}, $\theta_i = (\a_0 \meet \g^i_1) \join (\g^i_0 \meet \a_1) < \a_0 \meet \a_1$.

       Thus, the sublattice of $\con(\alg A(\b))$ generated by $\{\a_0, \a_1,  \b_0, \g_0^0, \g_1^0\}$ is isomorphic to $\alg K_{i+1}$.

     Case $\g^{i} <\g^{i+1}$ for all $i \in \mathbb{N}_0$: in this case we have an infinite chain of congruences $\{\g^j\}_{j \in \mathbb N_0}$ that clearly forms a $\alg K_{\infty}$. 
 \end{proof}

 Note that \cref{{thm:N5A(b)}} refines \cite[Theorem 1.2]{ABF} by providing a complete characterization of the sublattice of $\con(\alg A(\b))$ generated by $\{\a_0,$ $  \a_1,  \b_0, \g_0^0,$ $ \g_1^0\}$, removing all assumptions on $\alg N_5 \leq \con(\alg A)$. Furthermore, \cref{thm:N5A(b)} directly yields a characterization of the class of modular varieties in terms of lattice omission.

\begin{Cor}\label{cor:omisN5}
	The filter of omissions: $$\mathfrak{F}(\alg N_5) = \bigcap_{i \in \mathbb N} \mathfrak{F}(\alg M_i)\bigcap_{i \in \mathbb N} \mathfrak{F}(\alg K_i) \cap \mathfrak{F}(\alg K_{\infty}).$$
\end{Cor}

\begin{proof}
    The proof follows directly from \cref{thm:N5A(b)} observing that $\alg N_5$ is a sublattice of $\alg M_i$, $\alg K_i$ for all $i \in \mathbb{N}_0$, and of $\alg K_{\infty}$.
\end{proof}

\begin{Cor}\label{cor:N5char}
    Let $\vv V$ be a variety. Then the following conditions are equivalent:
    \begin{enumerate}
        \item $\vv V$ is congruence modular;
        \item $\vv V$ omits $\alg N_5$;
        \item $\vv V$ omits $\alg M_i$ and $\alg K_i$ for all $i \in \mathbb N$ and $\alg K_{\infty}$;
        \item $\vv V$ omits $\alg K$.
    \end{enumerate}
\end{Cor}

\begin{proof}
    The equivalence of $(1)$ and $(2)$ follows from \cite{Dedekind1900}, while the equivalence of $(2)$ and $(3)$ follows from \cref{cor:omisN5}. The equivalence of $(1)$ and $(4)$ follows from \cite[Theorem 6.100]{freese2022algebras}
\end{proof}

Note that the equivalences of \cref{cor:N5char} can also be proved without Freese's technique. To conclude the section, we provide examples of algebras satisfying the first two cases of \cref{thm:N5A(b)} for $i = 1$. These examples can be easily generalized to any $i \in \mathbb{N}_0$ as well as to the infinite case, thereby providing algebras satisfying each of the cases in \cref{thm:N5A(b)}.

  \begin{Ex}\label{ex:N5}
     We present examples of algebras that serve as instances for each of the first two cases in  \cref{thm:N5A(b)}, for $i =1$.
     \begin{enumerate}
         \item Let $\alg A$ be the algebra whose base set is $A = \{1,2,3,4,5,6, a,b,c,$ $d\}$ with only projections as operations. Then the partitions $\beta = \mid 1, 2 \mid 3, 4 \mid 5, 6 \mid a, c \mid b, d \mid$, $\gamma^0 = \mid 2, 3 \mid 4, 5 \mid a, b \mid c \mid d \mid$, $\gamma^1 = \mid 2, 3 \mid 4, 5 \mid a, b \mid c , d\mid$, and $\a = \mid 2, 3 \mid 4, 5 \mid 1, 6 \mid a, b \mid c, d\mid$   forms $\mathbf{N}_5 \leq \con(\alg A)$ labeled as in \cref{fig:N5D1D2} with $[\g, \a] = \{\g^0, \g^1, \a\}$ and, by \cref{thm:N5A(b)}, the sublattice of $\con(\alg A(\b))$ generated by $\{\a_0, \a_1,  \b_0, \g_0^0, \g_1^0\}$ is isomorphic to $\alg M_2$;

         \item let $\alg A$ be the algebra whose base set is $A = \{u_0,u_1,v_0,v_1,v_2, t_0,$ $t_1,t_2, t_3\}$ with only projections as operations. Then the partitions $\beta = \mid t_0,v_0,u_0 \mid t_1,v_2,t_2 \mid t_3,v_1, u_1 \mid$, $\gamma^0 = \mid t_0,t_1 \mid t_2, t_3 \mid$, $\gamma^1 = \mid v_0, v_1, v_2\mid t_0,t_1 \mid t_2, t_3 \mid$, and $\a = \mid u_0, u_1\mid v_0, v_1,v_2\mid t_0,t_1 \mid t_2, t_3 \mid$ forms $\mathbf{N}_5 \leq \con(\alg A)$ labeled as in \cref{fig:N5D1D2} with $[\g, \a] = \{\g^0, \g^1, \a\}$ and, by \cref{thm:N5A(b)}, the sublattice of $\con(\alg A(\b))$ generated by $\{\a_0, \a_1,  \b_0, \g_0^0, \g_1^0\}$ is isomorphic to $\alg K_2$.
     \end{enumerate}
 \end{Ex}

To construct an example of an algebra satisfying (2) of \cref{thm:N5A(b)}, the key is to define a chain $\g^0 < \g^1 < \g^2 = \a$ satisfying the hypothesis of \cref{thm:N5A(b)}. A diagram illustrating the relationships between the elements involved is provided in \cref{fig:exa} and can be fully generalized to construct an infinite chain, exhibiting exponential complexity in the diagram structure.

\begin{figure}[htbp]
		\begin{center}
			\begin{tikzcd}
				u_0 \arrow[dd,"\b",dash]\arrow[r,"\a", dash]&u_1\arrow[dd,"\b",dash]\\
				\\
				v_0\arrow[r,"\g^{1}", dash]&v_1
			\end{tikzcd}\\
            \begin{tikzcd}
				v_0 \arrow[dd,"\b",dash]\arrow[r,"\a", dash]&v_2\arrow[dd,"\b",dash]\arrow[r,"\a", dash]&v_2\arrow[dd,"\b",dash]\arrow[r,"\a", dash]&v_1\arrow[dd,"\b",dash]\\
				\\
				t_0\arrow[r,"\g^{0}", dash]&t_1&t_2\arrow[r,"\g^{0}", dash] &t_3.
			\end{tikzcd}
            \end{center}
        \caption{Graphical explanation of \cref{ex:N5} (2)}
        \label{fig:exa}
		\end{figure}

\section{An application of Freese's technique to $\alg D_1$}\label{sec:D1}

In this section, our goal is to characterize the omission of $\alg D_1$ with respect to the omission of $\alg D_{13}$ (\cref{fig:D13}). We will prove that the omissions of these lattices are equivalent, thereby providing another example of an application of Freese's technique producing a unique congruence lattice, similarly to the case of $\alg A(\g)$ for $\alg N_5$ (\cref{fig:N5D1D2}) and in contrast to what we have just seen in the previous section by considering $\alg A(\b)$. Our interest in any application of Freese's technique to $\alg D_1$ arises from the fact that the omission of this lattice characterizes Taylor varieties \cite{KearnesKiss2013}. Notably, the class of Taylor varieties forms a strong Mal'cev class \cite{Ols17}, and this class has a deep connection with the Feder-Vardi conjecture, which was independently proven in \cite{bulatov-dichotomy} and \cite{zhuk-dichotomy-short}. We demonstrate, using Freese's technique, that a variety is Taylor if and only if it omits $\alg D_{13}$, \cref{fig:D13}, a similar but different result with respect to \cite[Theorem 4.23]{KearnesKiss2013}.

\begin{figure}[htbp]
		\begin{center}
			\begin{tikzpicture}[scale=0.65]				
				\draw (-8,6) --(0,0) -- (8,6) -- (0,12) -- (-8,6);
                \draw (0,0) -- (-2,3) -- (-8,6);
                \draw (0,0) -- (2,3) -- (8,6);
                \draw (-2,3)-- (4,6) --(4,3) -- (0,6) -- (-4,3) -- (-4,6) -- (2,3);
                \draw (-4,6) -- (4,9) -- (0,6);
                \draw (4,6) -- (-4,9) -- (0,6);

                \draw[fill] (0,0) circle [radius=0.1];
				\node [left] at (0,0) {\footnotesize $0_{\alg A(\b)}$};
    
                \draw[fill] (-2,3) circle [radius=0.1];
				\node [left] at (-2,3) {\footnotesize $\a_0 \wedge \a_1$};
                 \draw[fill] (2,3) circle [radius=0.1];
				\node [right] at (2.05,3) {\footnotesize $\g_0 \wedge \g_1$};

                \draw[fill] (-4,3) circle [radius=0.1];
				\node [left] at (-4.1,2.95) {\footnotesize $\eta_0$};

                \draw[fill] (4,3) circle [radius=0.1];
				\node [right] at (4.1,2.95) {\footnotesize $\eta_1$};

                \draw[fill] (0,6) circle [radius=0.1];
				\node [right] at (0.1,6) {\footnotesize $\b_0$};

                \draw[fill] (-4,6) circle [radius=0.1];
				\node [left] at (-4.1,6) {\footnotesize $\g_0$};

                \draw[fill] (4,6) circle [radius=0.1];
				\node [right] at (4.1,6) {\footnotesize $\a_1$};

                 \draw[fill] (-8,6) circle [radius=0.1];
				\node [left] at (-8.1,6) {\footnotesize $\a_0$};

                \draw[fill] (8,6) circle [radius=0.1];
				\node [right] at (8.1,6) {\footnotesize $\g_1$};

                \draw[fill] (-4,9) circle [radius=0.1];
				\node [left] at (-4.1,9) {\footnotesize $\mu_0$};

                \draw[fill] (4,9) circle [radius=0.1];
				\node [right] at (4.1,9) {\footnotesize $\d_0$};

                \draw[fill] (0,12) circle [radius=0.1];
				\node [left] at (0,12.1) {\footnotesize $1_0 = 1_1$};

			\end{tikzpicture}
		\end{center}
		\caption{$\alg D_{13}$}\label{fig:D13}
	\end{figure}

 \begin{Thm}\label{thm:K13D1}
      Let $\alg A$ be an algebra and let $\con(\alg A)$ be its congruence lattice with a sublattice isomorphic to $\alg D_1$, labeled as in \cref{fig:N5D1D2}. Then, the sublattice of $\con(\alg A(\b))$ generated by $\{\a_0, \a_1, \g_0, \g_1,\b_0\}$ is isomorphic to $\alg D_{13}$ (\cref{fig:D13}).
 \end{Thm}

 \begin{proof}
 Again, as observed in \cref{sec:FreeseTech}, we may assume without loss of generality that the bottom element of $\alg D_1$ coincides with the bottom element of $\con(\alg A)$.
     		Then, we can observe that $\op{Con}(\alg A(\g)/\eta_i)$ is isomorphic to $\con(\alg A)$ through the natural isomorphism $f(\lambda) = \lambda_i$. Thus, the joins and meets in the part above $\eta_i$ in Figure \ref{fig:D13} are all respected, for $i \in [1]$. For the remaining part of the lattice, we can observe that $\a_0 \wedge \a_1, \g_0 \wedge \g_1 > 0_{\alg A(\b)}$ since $\a,\g >0_{\alg A}$. Furthermore, $\a_0 \wedge \a_1 \wedge\g_0 \wedge \g_1 = \a_0 \wedge \g_0 \wedge \a_1 \wedge \g_1 = \eta_0 \wedge \eta_1 = 0_{\alg A(\b)}$. By \cref{lem:etameet}, we have that $\eta_0 \wedge \a_0 \wedge \a_1 = \eta_0 \wedge \a_1 =  0_{\alg A}$ and similarly $\eta_0 \wedge \g_0 \wedge \g_1 =  0_{\alg A}$. Moreover, $\theta_0 \wedge \a_0 \wedge \a_1 = \eta_0 \wedge \a_0 \wedge \a_1 = 0_{\alg A(\b)}$, for all $\theta \in \{\eta, \beta, \g,\d\}$ while $\theta_0 \wedge \a_0 \wedge \a_1 = \a_0 \wedge \a_1$, for $\theta \in \{\a, \mu,  1_{\alg A(\b)}\}$. In a similar way, we can also prove that the meet between $\a_0 \wedge \a_1$ and the elements of the principal filter generated by $\eta_1$ and between $\g_0 \wedge \g_1$ and the elements of the principal filters generated by $\eta_0$ and $\eta_1$ are also correctly depicted. Regarding the joins, $\theta_0 \vee (\a_0 \wedge \a_1) = \theta_0 \vee \eta_0 \vee (\a_0 \wedge \a_1) = \theta_0 \vee \alpha_0 = (\theta \vee \a)_0$, for all $\theta \in \con(\alg A)$. Hence, the joins between $\a_0 \wedge \a_1$ and the elements of the filter generated by $\eta_0$ are correctly depicted. Similarly, the remaining joins can also be shown to be correctly depicted in \cref{fig:D13}.
 \end{proof}

      \begin{Thm}\label{TaylorThmorig}
		Let $\vv{V}$ be a variety. Then the following are equivalent:	
		\begin{enumerate}		
			\item $\vv{V}$ is Taylor, i. e., $\vv V$ satisfies a non-trivial identity;
			\item $\vv{V}$ satisfies the term equations:
			\begin{equation}\label{Olsak eq}
				t(x,y,y,y,x,x) \approx  t(y,x,y,x,y,x) \approx 	t(y,y, x,x,x,y);
			\end{equation}
			for some idempotent term $t$;
			\item $\vv{V}$ satisfies the following congruence equation in the variables $\{\alpha_1, \dots, \alpha_6,$ $ \beta_1, \dots, \beta_6\}$:
			\begin{equation}\label{Tayeq}
				\bigwedge^{6}_{i = 1}(\a_{i} \circ \b_{i}) \leq (\bigvee^{6}_{i=1}\a_{i} \wedge \bigwedge^{2}_{i=1} (\g \vee \theta_{i})) \vee (\bigvee^{6}_{i=1}\b_{i} \wedge \bigwedge^{2}_{i=1} (\g \vee \theta_{i}));
			\end{equation}
			where $\gamma =  \bigwedge^{6}_{i = 1} (\a_{i} \vee \beta_i)$ and 
			\begin{equation*}
				\theta_{i}= (\bigvee_{j \in L_{i}} \a_{j} \vee \bigvee_{j \in L'_{i}} \b_{j}) \wedge (\bigvee_{j \in R_{i}} \a_{j} \vee \bigvee_{j \in R'_{i}} \b_{j})
			\end{equation*}
			with 
			\begin{align*}
				&L_1 = L_2 =\{1,5,6\}  &L'_{1} = L'_2 = \{ 2,3,4 \}
				\\&R_{1} = \{ 2,4,6\} &R'_{1} = \{1,3,5 \}
				\\&R_{2} = \{ 3,4,5\} &R'_{2} = \{1,2,6 \};
			\end{align*}
            \item $\vv{V}$ satisfies the following congruence equation:
			\begin{equation*}\label{Tayeqweak}
				\bigwedge^{6}_{i = 1}(\a_{i} \circ \b_{i}) \leq (\bigvee^{6}_{i=1}\a_{i} \wedge \bigwedge^{2}_{i=1} (\tau \vee \theta_{i})) \vee (\bigvee^{6}_{i=1}\b_{i} \wedge \bigwedge^{2}_{i=1} (\tau \vee \theta_{i}));
			\end{equation*}
			where 
			\begin{equation*}
				\theta_{i}=(\bigvee_{j \in L_{i}} \a_{j} \vee \bigvee_{j \in L'_{i}} \b_{j}) \wedge ( [\tau,\tau] \circ \bigvee_{j \in R_{i}} \a_{j} \vee \bigvee_{j \in R'_{i}} \b_{j})
			\end{equation*}
			where $[,]$ is the TC-commutator and with $\tau = \bigwedge^{6}_{i = 1}(\a_{i} \vee\b_{i})$ and
			\begin{align*}
				&L_1 = L_2 =\{1,5,6\}  &L'_{1} = L'_2 = \{ 2,3,4 \}
				\\&R_{1} = \{ 2,4,6\} &R'_{1} = \{1,3,5 \}
				\\&R_{2} = \{ 3,4,5\} &R'_{2} = \{1,2,6 \}
			\end{align*}
            \item $\vv{V}(\alg A )\in \mathfrak{F}(\alg D_1)$;
            \item $\vv{V}(\alg A)\in \mathfrak{F}(\alg D_{13})$.         
		\end{enumerate}
	\end{Thm}

    \begin{proof}
        The equivalence of $(1)$ and $(2)$ follows from \cite{Ols17},  $(1),(5)$ are equivalent by \cite[Theorem $4.23$]{KearnesKiss2013}. The equivalence of $(1), (3)$, and $(4)$ follows from \cite{Fioravanti20241273}. For the last equivalence, we can observe that, as a consequence of \cref{thm:K13D1}, the omission of $\alg D_{13}$ implies the omission of $\alg D_1$, while the reverse is an immediate consequence of the fact that $\alg D_1 \leq \alg D_{13}$, and thus the claim holds.
    \end{proof}

    The preceding theorem illustrates how Freese's technique, originally applied to $\alg N_5$, can be extended to obtain a uniquely determined congruence lattice of a congruence viewed as a subalgebra of the square of the algebra, even for more complex lattices. Observe that there are several other conditions equivalent to those in \cref{TaylorThmorig}; we have presented only some of them. Note that the equivalence of items $(5)$ and $(6)$ of the previous theorem can also be proven without using Freese's technique, noting that $\alg D_{13} \leq \alg D_2^2$.

\section{An application of Freese's technique to $\alg D_2$}\label{sec:D2}
	
		In this section, our aim is to apply Freese's technique to the lattice $\alg D_2$ (\cref{fig:N5D1D2}) to investigate the structure of the congruence lattice $\con(\alg A(\mu))$. In order to do so, we introduce two families of lattices $\{\alg S_i\}_{i \in \mathbb{N}}$ and $\{\alg S_i^*\}_{i \in \mathbb{N}}$ (\cref{S_i}) along with the lattice $\alg S_{\infty}$ (\cref{Sinf}),  where $\theta_i = (\d_0^i \wedge \d_1^{i+1}) \vee (\d_0^{i+1} \wedge \d_1^{i})$, for all $i \in \mathbb N$. 
        For the whole section, we will use \cref{def:cong_chain} with $\delta^0(\g,\mu,\delta) = \delta$, $\delta^{i + 1}(\g,\mu,\delta) =  \Cg((\mu \circ \d^i \circ \mu) \cap \g)$.
	
	   We will prove that the omission of these two families of lattices, together with the lattice $\alg S_{\infty}$, is equivalent to the omission of $\alg D_2$. Furthermore, we will see that the presence of any of these lattices is intrinsically related to the behavior of the relational product of $\delta$ and $\mu$.

       \begin{figure}[htbp]
		\begin{center}
			\begin{tikzpicture}[scale=0.75]				
				\draw (-3,2) --(0,0) -- (3,2);
				\draw (-3,3) --(0,1) -- (3,3);
				\draw (-3,5) --(0,3) -- (3,5);
				\draw (-3,12) --(0,14) -- (3,12);
				
				\draw (-3,7) --(0,5) -- (3,7);
				\draw (-3,8) --(0,6) -- (3,8);
				\draw (-3,10) --(0,8) -- (3,10);
				\draw (-3,12) --(0,10) -- (3,12);
				\draw (-3,2) --(0,11) -- (3,2);
				\draw (-3,3) --(0,12) -- (3,3);
				\draw (0,11) --(-1.5,12) -- (0,14);
				
				\draw (0,0) --(0,1);
				\draw (0,2) --(0,3);
				\draw (3,2) --(3,5);
				\draw (-3,2) --(-3,5);
				\draw (3,5) --(3,7);
				\draw (-3,5) --(-3,7);
				\draw (0,4) --(0,5);
				\draw (0,7) --(0,8);
				\draw (0,8) --(0,10);
				\draw (0,11) --(0,14);
				\draw (3,8) --(3,12);
				\draw (-3,8) --(-3,12);
				
				\draw [dashed](-3,7) --(-3,8);
				\draw [dashed](3,7) --(3,8);
				\draw [dashed](0,5) --(0,6);
				
				\draw(-0.75,1.5) -- (0,2) -- (0.75,1.5);
				\draw(-0.75,3.5) -- (0,4) -- (0.75,3.5);
				\draw(-0.75,6.5) -- (0,7) -- (0.75,6.5);

				\draw[fill] (0,0) circle [radius=0.05];
				\node [left] at (0,0) {\footnotesize $0_{\alg A(\mu)}$};
				\draw[fill] (0,1) circle [radius=0.05];
				\node [right] at (0,1) {\tiny $\d^0_0 \wedge \d_0^1$};
				\draw[fill] (0,2) circle [radius=0.05];
				\node [right] at (0,2) {\tiny $\theta_0$};
				\draw[fill] (0,3) circle [radius=0.05];
				\node [right] at (0,3) {\tiny $\d^1_0 \wedge \d_1^1$};
				\draw[fill] (0,4) circle [radius=0.05];
				\node [right] at (0,4) {\tiny $\theta_1$};
				\draw[fill] (0,5) circle [radius=0.05];
				\node [right] at (0,5) {\tiny $\d^{2}_0 \wedge \d_1^{2}$};
				\node [right] at (0,6) {\tiny $\d^{i-1}_0 \wedge \d_1^{i-1}$};
				\draw[fill] (0,6) circle [radius=0.05];
				\node [right] at (0,7) {\tiny $\theta_{i-1}$};
				\draw[fill] (0,7) circle [radius=0.05];
				\node [right] at (0,8) {\tiny $\d^{i}_0 \wedge \d_1^{i}$};
				\draw[fill] (0,8) circle [radius=0.05];
				\node [right] at (0,10) {\tiny $\g_0 \wedge \g_1$};
				\draw[fill] (0,10) circle [radius=0.05];
				\node [right] at (0,12) {\tiny $\b_0 = \b_1$};
				\draw[fill] (0,12) circle [radius=0.05];
				\node [left] at (-1.4,12) {\tiny $\a_0 = \a_1$};
				\draw[fill] (-1.5,12) circle [radius=0.05];
				\node [right] at (0,14) {\tiny $1_0 = 1_1$};
				\draw[fill] (0,14) circle [radius=0.05];
				 \node [right] at (0,11) {\tiny $\mu_0 = \mu_1$};
				 \draw[fill] (0,11) circle [radius=0.05];

				\draw[fill] (-0.75,1.5) circle [radius=0.05];
				\node [left] at (-0.85,1.5) {\tiny $\d^{0}_0 \wedge \d_1^{1}$};
				\draw[fill] (0.75,1.5) circle [radius=0.05];
				\node [right] at (0.85,1.5) {\tiny $\d^{1}_0 \wedge \d_1^{0}$};
				\draw[fill] (-0.75,3.5) circle [radius=0.05];
				\node [left] at (-0.85,3.5) {\tiny $\d^{1}_0 \wedge \d_1^{2}$};
				\draw[fill] (0.75,3.5) circle [radius=0.05];
				\node [right] at (0.85,3.5) {\tiny $\d^{2}_0 \wedge \d_1^{1}$};
				\draw[fill] (-0.75,6.5) circle [radius=0.05];
				\node [left] at (-0.85,6.5) {\tiny $\d^{i-1}_0 \wedge \d_1^{i}$};
				\draw[fill] (0.75,6.5) circle [radius=0.05];
				\node [right] at (0.85,6.5) {\tiny $\d^{i}_0 \wedge \d_1^{i-1}$};
				
				\draw[fill] (3,2) circle [radius=0.05];
				\node [right] at (3,2) {\tiny $\eta_1$};
				\draw[fill] (3,3) circle [radius=0.05];
				\node [right] at (3,2.8) {\tiny $\d^0_1$};
				\draw[fill] (3,5) circle [radius=0.05];
				\node [right] at (3,5) {\tiny $\d^1_1$};
				\draw[fill] (3,7) circle [radius=0.05];
				\node [right] at (3,7) {\tiny $\d^{2}_1$};
				\draw[fill] (3,8) circle [radius=0.05];
				\node [right] at (3,8) {\tiny $\d^{i-1}_1$};
				\draw[fill] (3,10) circle [radius=0.05];
				\node [right] at (3,10) {\tiny $\d^{i}_1$};
				\draw[fill] (3,12) circle [radius=0.05];
				\node [right] at (3,12) {\tiny $\g_1$};

				\draw[fill] (-3,2) circle [radius=0.05];
				\node [left] at (-3,2) {\tiny $\eta_0$};
				\draw[fill] (-3,3) circle [radius=0.05];
				\node [left] at (-3,3) {\tiny $\delta^0_0$};
				\draw[fill] (-3,5) circle [radius=0.05];
				\node [left] at (-3,5) {\tiny $\d^1_0$};
				\draw[fill] (-3,7) circle [radius=0.05];
				\node [left] at (-3,7) {\tiny $\d^{2}_0$};
				\draw[fill] (-3,8) circle [radius=0.05];
				\node [left] at (-3,8) {\tiny $\d^{i-1}_0$};
				\draw[fill] (-3,10) circle [radius=0.05];
				\node [left] at (-3,10) {\tiny $\d^{i}_0$};
				\draw[fill] (-3,12) circle [radius=0.05];
				\node [left] at (-3,12) {\tiny $\g_0$};


				\draw (5,2) --(8,0) -- (11,2);
				\draw (5,3) --(8,1) -- (11,3);
				\draw (5,5) --(8,3) -- (11,5);
				\draw (5,12) --(8,14) -- (11,12);
				
				\draw (5,7) --(8,5) -- (11,7);
				\draw (5,8) --(8,6) -- (11,8);
				\draw (5,10) --(8,8) -- (11,10);
				\draw (5,12) --(8,10) -- (11,12);
				\draw (5,2) --(8,11) -- (11,2);
				\draw (5,3) --(8,12) -- (11,3);
				\draw (8,11) --(6.5,12) -- (8,14);
				
				\draw (8,0) --(8,1);
				\draw (8,2) --(8,3);
				\draw (11,2) --(11,5);
				\draw (5,2) --(5,5);
				\draw (11,5) --(11,7);
				\draw (5,5) --(5,7);
				\draw (8,4) --(8,5);
				\draw (8,7) --(8,8);
				\draw (8,9) --(8,10);
				\draw (8,11) --(8,14);
				\draw (11,8) --(11,12);
				\draw (5,8) --(5,12);
				
				\draw [dashed](5,7) --(5,8);
				\draw [dashed](11,7) --(11,8);
				\draw [dashed](8,5) --(8,6);
				
				\draw(7.25,1.5) -- (8,2) -- (8.75,1.5);
				\draw(7.25,3.5) -- (8,4) -- (8.75,3.5);
				\draw(7.25,6.5) -- (8,7) -- (8.75,6.5);
				\draw(7.25,8.5) -- (8,9) -- (8.75,8.5);

				\draw[fill] (8,0) circle [radius=0.05];
				\node [left] at (8,0) {\footnotesize $0_{\alg A(\mu)}$};
				\draw[fill] (8,1) circle [radius=0.05];
				\node [right] at (8,1) {\tiny $\d^0_0 \wedge \d_0^1$};
				\draw[fill] (8,2) circle [radius=0.05];
				\node [right] at (8,2) {\tiny $\theta_0$};
				\draw[fill] (8,3) circle [radius=0.05];
				\node [right] at (8,3) {\tiny $\d^1_0 \wedge \d_1^1$};
				\draw[fill] (8,4) circle [radius=0.05];
				\node [right] at (8,4) {\tiny $\theta_1$};
				\draw[fill] (8,5) circle [radius=0.05];
				\node [right] at (8,5) {\tiny $\d^{2}_0 \wedge \d_1^{2}$};
				\node [right] at (8,6) {\tiny $\d^{i-1}_0 \wedge \d_1^{i-1}$};
				\draw[fill] (8,6) circle [radius=0.05];
				\node [right] at (8,7) {\tiny $\theta_{i-1}$};
				\draw[fill] (8,7) circle [radius=0.05];
				\node [right] at (8,8) {\tiny $\d^{i}_0 \wedge \d_1^{i}$};
				\draw[fill] (8,9) circle [radius=0.05];
				\node [right] at (8,9) {\tiny $\theta_i$};
				\draw[fill] (8,8) circle [radius=0.05];
				\node [right] at (8,12) {\tiny $\b_0 \wedge \b_1$};
				\draw[fill] (8,12) circle [radius=0.05];
				\node [right] at (8,10) {\tiny $\g_0 \wedge \g_1$};
				\draw[fill] (8,10) circle [radius=0.05];
				\node [left] at (6.6,12) {\tiny $\a_0 = \a_1$};
				\draw[fill] (6.5,12) circle [radius=0.05];
				\node [right] at (8,14) {\tiny $1_0 = 1_1$};
				\draw[fill] (8,14) circle [radius=0.05];
				\node [right] at (8,11) {\tiny $\mu_0 = \mu_1$};
				\draw[fill] (8,11) circle [radius=0.05];

				\draw[fill] (7.25,1.5) circle [radius=0.05];
				\node [left] at (7.15,1.5) {\tiny $\d^{0}_0 \wedge \d_1^{1}$};
				\draw[fill] (8.75,1.5) circle [radius=0.05];
				\node [right] at (8.85,1.5) {\tiny $\d^{1}_0 \wedge \d_1^{0}$};
				\draw[fill] (7.25,3.5) circle [radius=0.05];
				\node [left] at (7.15,3.5) {\tiny $\d^{1}_0 \wedge \d_1^{2}$};
				\draw[fill] (8.75,3.5) circle [radius=0.05];
				\node [right] at (8.85,3.5) {\tiny $\d^{2}_0 \wedge \d_1^{1}$};
				\draw[fill] (7.25,6.5) circle [radius=0.05];
				\node [left] at (7.15,6.5) {\tiny $\d^{i-1}_0 \wedge \d_1^{i}$};
				\draw[fill] (8.75,6.5) circle [radius=0.05];
				\node [right] at (8.85,6.5) {\tiny $\d^{i}_0 \wedge \d_1^{i-1}$};
				\draw[fill] (7.25,8.5) circle [radius=0.05];
				\node [left] at (7.15,8.5) {\tiny $\d^{i}_0 \wedge \g_1$};
				\draw[fill] (8.75,8.5) circle [radius=0.05];
				\node [right] at (8.85,8.5) {\tiny $\g_0 \wedge \d_1^{i}$};
				
				\draw[fill] (11,2) circle [radius=0.05];
				\node [right] at (11,2) {\tiny $\eta_1$};
				\draw[fill] (11,3) circle [radius=0.05];
				\node [right] at (11,3) {\tiny $\d^0_1$};
				\draw[fill] (11,5) circle [radius=0.05];
				\node [right] at (11,5) {\tiny $\d^1_1$};
				\draw[fill] (11,7) circle [radius=0.05];
				\node [right] at (11,7) {\tiny $\d^{2}_1$};
				\draw[fill] (11,8) circle [radius=0.05];
				\node [right] at (11,8) {\tiny $\d^{i-1}_1$};
				\draw[fill] (11,10) circle [radius=0.05];
				\node [right] at (11,10) {\tiny $\d^{i}_1$};
				\draw[fill] (11,12) circle [radius=0.05];
				\node [right] at (11,12) {\tiny $\g_1$};

				\draw[fill] (5,2) circle [radius=0.05];
				\node [left] at (5,2) {\tiny $\eta_0$};
				\draw[fill] (5,3) circle [radius=0.05];
				\node [left] at (5,3.2) {\tiny $\delta^0_0$};
				\draw[fill] (5,5) circle [radius=0.05];
				\node [left] at (5,5) {\tiny $\d^1_0$};
				\draw[fill] (5,7) circle [radius=0.05];
				\node [left] at (5,7) {\tiny $\d^{2}_0$};
				\draw[fill] (5,8) circle [radius=0.05];
				\node [left] at (5,8) {\tiny $\d^{i-1}_0$};
				\draw[fill] (5,10) circle [radius=0.05];
				\node [left] at (5,10) {\tiny $\d^{i}_0$};
				\draw[fill] (5,12) circle [radius=0.05];
				\node [left] at (5,12) {\tiny $\g_0$};
				
			\end{tikzpicture}
		\end{center}
		\caption{$\alg S_{i+1}$ and $\alg S_{i+1}^*$}\label{S_i}
	\end{figure}
	
		\begin{figure}[htbp]
		\begin{center}
			\begin{tikzpicture}[scale=0.90]
				
				\draw (-3,2) --(0,0) -- (3,2);
				\draw (-3,3) --(0,1) -- (3,3);
				\draw (-3,5) --(0,3) -- (3,5);

				\draw (-3,7) --(0,5) -- (3,7);
				\draw (-3,8) --(0,6) -- (3,8);
				\draw (3,2) --(0,7)--(-3,2);
				\draw (3,3) --(0,8)--(-3,3);
				\draw (3,8)--(0,10)--(-3,8);
				\draw (0,7) --(-1.5,8)--(0,10);

				\draw (0,0) --(0,1);
				\draw (0,2) --(0,3);
				\draw (3,2) --(3,5);
				\draw (-3,2) --(-3,5);
				\draw (3,5) --(3,7);
				\draw (-3,5) --(-3,7);
				\draw (0,4) -- (0,5);
				\draw (0,7) --(0,10);

				\draw [dashed](-3,7) --(-3,8);
				\draw [dashed](3,7) --(3,8);
				\draw [dashed](0,5) --(0,6);
				
				\draw(-0.75,1.5) -- (0,2) -- (0.75,1.5);
				\draw(-0.75,3.5) -- (0,4) -- (0.75,3.5);

				\draw[fill] (0,0) circle [radius=0.05];
				\node [left] at (0,0) {\footnotesize $0_{\alg A(\mu)}$};
				\draw[fill] (0,1) circle [radius=0.05];
				\node [right] at (0,1) {\tiny $\d^0_0 \wedge \d_1^0$};
				\draw[fill] (0,2) circle [radius=0.05];
				\node [right] at (0,2) {\tiny $\theta_0$};
				\draw[fill] (0,3) circle [radius=0.05];
				\node [right] at (0,3) {\tiny $\d^1_0 \wedge \d_1^1$};
				\draw[fill] (0,4) circle [radius=0.05];
				\node [right] at (0,4) {\tiny $\theta_1$};
				\draw[fill] (0,5) circle [radius=0.05];
				\node [right] at (0,5) {\tiny $\d^{2}_0 \wedge \d_1^{2}$};
				\node [right] at (0,6) {\tiny $\g_0 \wedge \g_1$};
				\draw[fill] (0,6) circle [radius=0.05];
				\node [right] at (0,8) {\tiny $\b_0 = \b_1$};
				\draw[fill] (0,8) circle [radius=0.05];
				\node [right] at (0,7) {\tiny $\mu_0 = \mu_1$};
				\draw[fill] (0,7) circle [radius=0.05];
				\node [right] at (0,10) {\tiny $1_0 = 1_1$};
				\draw[fill] (0,10) circle [radius=0.05];

				\draw[fill] (-0.75,1.5) circle [radius=0.05];
				\node [left] at (-0.85,1.5) {\tiny $\d^{0}_0 \wedge \d_1^{1}$};
				\draw[fill] (0.75,1.5) circle [radius=0.05];
				\node [right] at (0.85,1.5) {\tiny $\d^{1}_0 \wedge \d_1^{0}$};
				\draw[fill] (-0.75,3.5) circle [radius=0.05];
				\node [left] at (-0.85,3.5) {\tiny $\d^{1}_0 \wedge \d_1^{2}$};
				\draw[fill] (0.75,3.5) circle [radius=0.05];
				\node [right] at (0.85,3.5) {\tiny $\d^{2}_0 \wedge \d_1^{1}$};

				\draw[fill] (3,2) circle [radius=0.05];
				\node [right] at (3,2) {\tiny $\eta_1$};
				\draw[fill] (3,3) circle [radius=0.05];
				\node [right] at (3,3) {\tiny $\d^0_1 $};
				\draw[fill] (3,5) circle [radius=0.05];
				\node [right] at (3,5) {\tiny $\d^1_1$};
				\draw[fill] (3,7) circle [radius=0.05];
				\node [right] at (3,7) {\tiny $\d^{2}_1$};
				\draw[fill] (3,8) circle [radius=0.05];
				\node [right] at (3,8) {\tiny $\g_1$};

				\draw[fill] (-3,2) circle [radius=0.05];
				\node [left] at (-3,2) {\tiny $\eta_0$};
				\draw[fill] (-3,3) circle [radius=0.05];
				\node [left] at (-3,3) {\tiny $\d^0_0 $};
				\draw[fill] (-3,5) circle [radius=0.05];
				\node [left] at (-3,5) {\tiny $\d^1_0$};
				\draw[fill] (-3,7) circle [radius=0.05];
				\node [left] at (-3,7) {\tiny $\d^{2}_0$};
				\draw[fill] (-3,8) circle [radius=0.05];
				\node [left] at (-3,8) {\tiny $\g_0$};
				\draw[fill] (-1.5,8) circle [radius=0.05];
				\node [left] at (-1.4,8) {\tiny $\a_0 = \a_1$};

			\end{tikzpicture}
		\end{center}
		\caption{$\alg S_{\infty}$}\label{Sinf}
	\end{figure}
	
	\begin{Thm}\label{thm:D2}
		Let $\alg A$ be an algebra such that $\con(\alg A)$ has a sublattice as the $\alg D_2$ in Figure \ref{fig:N5D1D2}. Let us consider $\con(\alg A(\mu))$ and let $\{\delta^{i }(\g,\mu,\delta)\}_{i \in \mathbb{N}_0}$ be a chain of congruences as defined in \ref{def:cong_chain}. Then the sublattice of $\con(\alg A(\mu))$ generated by $\{\mu_0, \a_0,  \g_0, \g_1,$ $ \d_0^0, \d_1^0\}$ is isomorphic to:
    \begin{enumerate}
        \item $\alg S_{i+1}$ if and only if $\d^{i-1} <\d^{i+1} = \d^i <\g$;
        \item $\alg S^*_{i+1}$ if and only if $\d^{i} < \d^{i+1} =\g$;
        \item $\alg S_{\infty}$ if and only if $\d^{i} < \d^{i+1}$ for all $i \in \mathbb{N}_0$.
    \end{enumerate}
    \end{Thm}

     \begin{proof}
     The strategy will be the same as \cref{thm:N5A(b)}. Again, as observed in \cref{sec:FreeseTech}, we may assume without loss of generality that the bottom element of $\alg D_2$ coincides with the bottom element of $\con(\alg A)$.
     		Then,  we make some observations about the common structure shared by all the cases. We note that $\op{Con}(\alg A(\mu)/\eta_s)$ is isomorphic to $\con(\alg A)$ through the natural isomorphism $f(\lambda) = \lambda_s$, for $s \in \{0,1\}$. Furthermore, if $\d^j < \d^{j+1}$, then also $\d^j_0 \wedge \d^j_1  < \d^{j+1}_0 \wedge \d^{j+1}_0$, for $j \in \mathbb N_0$. Moreover, by \cref{lem:skew} and \cref{lem:gen_gi}, whenever $\d^j <\d^{j+1}$, we have that 
     \[
     0 <\d^j_0 \wedge \d^{j+1}_1, \d^{j+1}_0 \wedge \d^{j}_1 < \d^{j+1}_0 \wedge \d^{j+1}_1.
     \]

   Hence, by \cref{lem:gen_gi}, we can observe that

     \[
     \eta_0 \vee (\g_0 \wedge \d_1^j) = \eta_0 \vee \d^j_0 \vee (\g_0 \wedge \d_1^j) = \d^{j+1}_0.
     \] 

     for all $j \in \mathbb{N}_0$. Moreover, $\d_0^{j+1} \wedge \d^j_1 = \g_0 \wedge \d^j_1$, for all $j \in \mathbb{N}_0$, follows from the \cref{def:cong_chain} since $((a,b),(x_1,x_2)) \in \g_0 \wedge \d^j_1$ implies $(a,b) \in \d^{j+1}$ and $(x_1,x_2) \in \d^j$. The other joins and meets follow from the order relations between the factor congruences.

     What is left to prove is basically which of the congruences in $\alg S_{\infty}$ is collapsing. By \cref{lemma1}, $\theta_j = (\d^{j+1}_0 \meet \d^j_1) \join (\d_0^j \meet \d^{j+1}_1) < \d^{j+1}_0 \meet \d^{j+1}_1$, for all $j \in \mathbb N_0$ with $\d^{i} <\d^{i+1}$.

     Case $\d^{i-1} <\d^{i+1} = \d^i <\g$: then, by \cref{lem:skew}, $\d^i_0 \wedge \g_1 = \d^i_0 \wedge \d^i_1$, the chain $\{\d^j\}_{j \in \mathbb N_0}$ collapses for $j \geq i$  and thus the sublattice of $\con(\alg A(\mu))$ generated by $\{\mu_0, \a_0,  \g_0, \g_1,$ $ \d_0^0, \d_1^0\}$ is isomorphic to $\alg S_{i+1}$.

     Case $\d^{i} <\d^{i+1} =\g$: then, by \cref{lem:skew}, $\d^i_0 \wedge \g_1 > \d^i_0 \wedge \d^i_1$, $\d_0^i \vee(\d^i_0 \wedge \g_1) = \g_0$, and the chain $\{\d^j\}_{j \in \mathbb N_0}$ collapses for $j > i$. Furthermore, by \cref{lemma1}, $\theta_i = (\g_0 \meet \d^i_1) \join (\d^i_0 \meet \g_1) < \g_0 \meet \g_1$.

     Thus, the sublattice of $\con(\alg A(\mu))$ generated by$\{\mu_0, \a_0,  \g_0, \g_1,$ $ \d_0^0, \d_1^0\}$ is isomorphic to $\alg S^*_{i+1}$.

     Case $\d^{i} <\d^{i+1}$ for all $i \in \mathbb{N}_0$: in this case, we have an infinite chain of congruences $\{\g^j\}_{j \in \mathbb N_0}$ which clearly forms a $\alg S_{\infty}$.
 \end{proof}

Note that by applying the same strategy used to construct \cref{ex:N5}, we can similarly generate algebras that satisfy each of the cases in \cref{thm:D2}, thereby demonstrating their non-emptiness.

 \begin{Cor}\label{cor:filD2}
	The filter of omissions $\mathfrak{F}(\alg D_2) = \bigcap_{i \in \mathbb N} \mathfrak{F}(\alg S_i)\bigcap_{i \in \mathbb N} \mathfrak{F}(\alg S^*_i) \cap \mathfrak{F}(\alg S_{\infty})$
\end{Cor}

\begin{proof}
    The proof follows from \cref{thm:D2} observing that $\alg D_2$ is a sublattice of both $\alg S^*_i$ and $\alg S_i$, for all $i \in \mathbb{N}$, and of $\alg S_{\infty}$.
\end{proof}

\begin{Thm}\label{thm:d2char}
    Let $\vv V$ be a variety. Then the following conditions are equivalent:
    \begin{enumerate}
        \item $\vv V$ satisfies a non-trivial congruence identity;
        \item $\vv V$ omits $\alg D_2$;
        \item $\vv V$ omits $\alg S_i$ and $\alg S_i^*$ for all $i \in \mathbb N$ and $\alg S_{\infty}$.
        
    \end{enumerate}
\end{Thm}

\begin{proof}
    The equivalence of $(1)$ and $(2)$ follows from \cite[Theorem 8.11] {KearnesKiss2013}, while the equivalence of $(2)$ and $(3)$ follows from \cref{cor:filD2}.
\end{proof}

	\section{Conclusions and problems}\label{sec:Conclusions}
	
    In this paper, we applied Freese’s technique to derive characterizations of congruence modular varieties, Taylor varieties, and varieties satisfying a non-trivial congruence identity through lattice omission, further demonstrating the effectiveness of this methodology.

    Regarding the future development of this work, we have observed that, depending on certain properties of the congruence sublattice considered, Freese's technique uniquely identifies a lattice whose omission is equivalent to that of the original (e.g. \cref{thm:K13D1}) or allows for multiple possibilities (e.g. \cref{thm:D2}). Understanding the precise conditions under which this technique yields a unique lattice would provide deeper insight into its underlying mechanisms.

    Additionally, several open questions arise from the concept of the filter of omission. As far as we know, there exist only four non-trivial Mal’cev classes that coincide with the omission class of some finite subdirectly irreducible lattice  $\alg L$:
	\begin{enumerate}
		\item the class $T$ of all Taylor varieties, i.e., satisfying a nontrivial idempotent Mal'cev condition (\cite{KearnesKiss2013}, Theorem 4.23);
		\item the class $L$ of all varieties $\vv V$ such that $\Con(\vv V)$ is a proper subvariety of the variety  of lattices (\cite{KearnesKiss2013}, Theorem 8.11), i.e., the varieties omitting $\alg D_2$;
		\item the class $M$ of all congruence modular varieties (really from \cite{Dedekind1900});
		\item the class $SD_\meet$ of all congruence meet semidistributive varieties (\cite{KearnesKiss2013}, Theorem 8.1).
	\end{enumerate}
	
	It is natural to ask whether these are the only Mal’cev classes with this property. A more challenging question is whether necessary and/or sufficient conditions can be established for a finite subdirectly irreducible lattice 
$\alg L$ such that the omission class of $\alg L$ forms a Mal’cev class.

As a future direction, we aim to apply Freese’s technique to other lattices whose omissions correspond to Mal’cev conditions. As mentioned in the Introduction, the problem of determining which lattices may or may not appear as congruence sublattices has a long history, dating back to Dedekind’s work \cite{Dedekind1900}. Despite the longevity of these questions, only a few examples of lattices have been identified whose omission characterizes significant Mal’cev conditions. Applying Freese’s technique to additional lattices may provide new insights into this field, potentially leading to novel characterizations of known Mal’cev conditions, as shown in the article, or the discovery of new interesting ones.

\subsection*{Ethical approval}
Not applicable. 

\subsection*{Competing interests} 
Not applicable. 

\subsection*{Authors' contributions} 
Not applicable. 

\subsection*{Availability of data and materials}
Not applicable.

\bibliographystyle{plain}
\bibliography{bibliografia.bib}

\end{document}